\documentclass[reqno]{amsart}

\newtheorem{Theorem}{Theorem}
\newtheorem{Proposition}[Theorem]{Proposition}
\newtheorem{Corollary}[Theorem]{Corollary}
\numberwithin{Theorem}{section}
\numberwithin{equation}{section}

\allowdisplaybreaks

\begin{document}

\title[Elementary derivations of bilateral series identities]
{Elementary derivations of identities for bilateral basic
hypergeometric series}
\author{Michael Schlosser}

\address{Department of Mathematics, The Ohio State University,
231 West 18th Avenue, Columbus, Ohio 43210, USA}
\email{mschloss@math.ohio-state.edu}
\date{November 1, 2000}
\thanks{2000 {\em Mathematics Subject Classification:} Primary 33D15.\\\indent
{\em Keywords and phrases:} bilateral basic hypergeometric series,
$q$-series, Bailey's $_6\psi_6$ summation,
Jackson's $_8\psi_8$ transformation, $_{10}\psi_{10}$ transformation,
Slater's transformations for well-poised $_{2r}\psi_{2r}$ series,
Slater's transformations for general $_{r}\psi_{r}$ series,
Chu--Gasper--Karlsson--Minton-type transformations.}

\dedicatory{{\rm\small Department of Mathematics, The Ohio State University,\\
231 West 18th Avenue, Columbus, Ohio 43210, USA\\
E-mail: \tt mschloss@math.ohio-state.edu}}

\begin{abstract}
We give elementary derivations of several classical and some new
summation and transformation formulae for bilateral basic
hypergeometric series. For purpose of motivation, we review
our previous simple proof (``A simple proof of Bailey's very-well-poised
$_6\psi_6$ summation", {\em Proc\@. Amer\@. Math\@. Soc\@.}, to appear)
of Bailey's
very-well-poised $_6\psi_6$ summation. Using a similar but different method,
we now give elementary derivations of some transformations for bilateral
basic hypergeometric series.
In particular, these include M.~Jackson's very-well-poised $_8\psi_8$
transformation, a very-well-poised
$_{10}\psi_{10}$ transformation, by induction,
Slater's general transformation for very-well-poised $_{2r}\psi_{2r}$ series,
and Slater's transformation for general $_{r}\psi_{r}$ series.
Finally, we derive some new transformations for bilateral
basic hypergeometric series of Chu--Gasper--Karlsson--Minton-type.
\end{abstract}

\maketitle

\section{Introduction}

The classical theories of hypergeometric series (cf.~\cite{slater})
and $q$-hypergeometric series (cf.~\cite{grhyp})
consist of many known summation and transformation formulae.
In fact, most identities for series involving 
($q$-)binomial coefficients can be uniformly written in terms of
($q$-)hypergeometric series.
Well-known examples are the binomial theorem, the Vandermonde summation,
and their ``$q$-analogues". There are numerous other summations,
and also transformations for ($q$-)hypergeometric series.
The $q$-hypergeometric series are usually called
basic hypergeometric series, where ``basic" refers to the base $q$.
The theory of basic hypergeometric series, which
contains hypergeometric series as special cases,
arose initially in combinatorics and classical analysis,
and interacts similarly with number theory, statistics, physics,
and representation theory of quantum Lie algebras,
see Andrews~\cite{qandrews}.

The theories of unilateral (or one-sided) hypergeometric and basic
hypergeometric series have quite a rich history, dating back to,
at least, Euler. Formulae for {\em bilateral} (basic) hypergeometric
series were not discovered until 1907 when Dougall~\cite{dougall},
using residue calculus, derived summations for the bilateral $_2H_2$ and
very-well-poised $_5H_5$ series.
Ramanujan~\cite{hardy} extended the $q$-binomial theorem by finding
a summation formula for the bilateral $_1\psi_1$ series.
Later, Bailey~\cite{bail66},\cite{bail22} carried out
systematical investigations on  bilateral basic hypergeometric
series. Further significant contributions were made by
Slater~\cite{slatertf},\cite{slater}, a student of Bailey.
See \cite{grhyp} and \cite{slater} for an excellent survey of
the above classical material.

Bailey's~\cite[Eq.~(4.7)]{bail66} very-well-poised $_6\psi_6$ summation
(cf.~\cite[Eq.~(5.3.1)]{grhyp}) is a very powerful identity,
as it stands at the top of the classical
hierarchy of summation formulae for bilateral series.
Some of the applications of the $_6\psi_6$ summation to partitions and
number theory are given in Andrews~\cite{andappl}.
Several proofs of Bailey's $_6\psi_6$ summation are already
known (see, e.g., Bailey~\cite{ bail66}, Slater and Lakin~\cite{slatlat},
Andrews~\cite{andappl}, Askey and Ismail~\cite{askmail}, and
Askey~\cite{askeyII}) which, unfortunately, are not entirely elementary.
Very recently, the author~\cite{schlelsum} found a new
simple proof of the very-well-poised
$_6\psi_6$ summation formula, directly from three applications of
Rogers'~\cite[p.~29, second eq.]{rogers} nonterminating
$_6\phi_5$ summation (cf.~\cite[Eq.~(2.7.1)]{grhyp}) and elementary
manipulations of series.

The method we used in \cite{schlelsum} extends that already used by
M.~Jackson~\cite[Sec.~4]{mjack} in her first elementary proof
of Ramanujan's $_1\psi_1$ summation formula~\cite{hardy}
(cf.~\cite[Eq.~(5.2.1)]{grhyp}). In \cite{schlelsum}, besides of giving
an elementary derivation of Bailey's very-well-poised $_6\psi_6$ summation,
we also gave an elementary derivation of Dougall's~\cite[Sec.~13]{dougall}
$_2H_2$ summation.

In this article, we apply a similar but different method
to derive several classical and some new
transformations for bilateral basic hypergeometric series.
In fact, here we make use of unilateral transformations and combine them
with bilateral series identities to deduce more complicated bilateral
series identities.
After recalling some standard notation for basic hypergeometric series
in Section~\ref{sec0}, we review, for purpose of motivation,
our~\cite{schlelsum} elementary proof of Bailey's
very-well-poised $_6\psi_6$ summation in Section~\ref{sec6}.
In Section~\ref{sec8},
we combine Bailey's~\cite{bail66} summation formula for a nonterminating
very-well-poised $_8\phi_7$ series (cf.~\cite[Eq.~(2.11.7)]{grhyp})
with Bailey's $_6\psi_6$ summation. As result we obtain a
transformation formula for a very-well-poised $_8\psi_8$ series
into a sum of two $_8\phi_7$ series.
This transformation is equivalent to a transformation given by
M.~Jackson~\cite[Eq.~(2.2)]{mjack1}.
In Section~\ref{sec10}, we apply our machinery to
deduce a transformation formula for a very-well-poised $_{10}\psi_{10}$
series into a sum of three $_{10}\phi_9$ series. This $_{10}\psi_{10}$
transformation is given implicitly by Slater~\cite{slatertf} and
explicitly by Gasper and Rahman~\cite[Eq.~(5.6.3)]{grhyp}.
We go even further in Section~\ref{sec2r}, where we prove
Slater's~\cite{slatertf} general transformation for very-well-poised
$_{2r}\psi_{2r}$ series by induction.
Similarly, in Section~\ref{secslgen}, we give an elementary
inductive derivation of Slater's~\cite{slatertf} general $_r\psi_r$
transformation. Instead of making use of Bailey's nonterminating
$_8\phi_7$ summation, we utilize the nonterminating $_3\phi_2$ summation here.
Finally, in Section~\ref{secc}, using Slater's
general transformations for bilateral basic hypergeometric series,
we give elementary derivations of transformations of
Chu--Gasper--Karlsson--Minton-type which seem to be new.
These identities generalize
some formulae recently found by Chu~\cite{chubs},\cite{chuwp}.

It is worth noting that most of the classical transformations for
bilateral basic hypergeometric series are proved in the literature
by specializing down from the very general transformations provided by
Slater~\cite{slatertf}. For instance, this is how the
very-well-poised $_8\psi_8$ and $_{10}\psi_{10}$ transformations
(see Equations~\eqref{88gl} and \eqref{1010gl}) are usually derived.
Slater derives her transformations in \cite{slatertf}
by using some general transformations
for basic hypergeometric series which Sears~\cite{searstf}
derived by manipulations of series.
Slater~\cite{slatertf2} also gives shorter proofs of her
transformations which use contour integrals of Barnes' type.
Already Watson~\cite{watson} had used such integrals
to derive transformations for basic hypergeometric series
of any order. In contrast, in this article we do not use contour integration.
After finding the bilateral transformations by elementary means,
we do sometimes appeal to analytic continuation to extend our results.
Our course of deriving the general
transformations~\eqref{2r2rgl} and \eqref{slgentfu}
is reminiscent of (but different from) Sears' analysis in \cite{searstf}.

The ideas in this article should open up new avenues in the theory of
multiple basic hypergeometric series. Whereas in the one-dimensional
theory it is possible to specialize down from general high level
identities, the extension of these to multiple basic hypergeometric series
are not yet known. We expect that our technique will allow to proceed from
lower level identities to systematically derive the upper level ones.
In particular, we plan to apply the methods of this article and of
\cite{schlelsum} in the setting of multiple
basic hypergeometric series associated to root systems, see e.g\@.
Milne~\cite{milne}, Gustafson~\cite{gus},
v.~Diejen~\cite{diejen}, and Schlosser~\cite{hypdet}.

Finally, we wish to gratefully acknowledge the helpful comments and
suggestions of George Andrews, Mourad Ismail, Christian Krattenthaler
and Stephen Milne.

\section{Background and notation}\label{sec0}

Here we recall some standard notation for $q$-series,
and basic hypergeometric series (cf. Gasper and Rahman~\cite{grhyp}).

Let $q$ be a complex number such that $0<|q|<1$. We define the
{\em $q$-shifted factorial} for all integers $k$ by
\begin{equation*}
(a;q)_{\infty}:=\prod_{j=0}^{\infty}(1-aq^j)\qquad\text{and}\qquad
(a;q)_k:=\frac{(a;q)_{\infty}}{(aq^k;q)_{\infty}}.
\end{equation*}
For brevity, we employ the usual notation
\begin{equation*}
(a_1,\ldots,a_m;q)_k\equiv (a_1;q)_k\dots(a_m;q)_k
\end{equation*}
where $k$ is an integer or infinity. Further, we utilize the notations
\begin{equation}\label{defhyp}
_r\phi_s\!\left[\begin{matrix}a_1,a_2,\dots,a_r\\
b_1,b_2,\dots,b_s\end{matrix};q,z\right]:=
\sum _{k=0} ^{\infty}\frac {(a_1,a_2,\dots,a_r;q)_k}
{(q,b_1,\dots,b_s;q)_k}\left((-1)^kq^{\binom k2}\right)^{1+s-r}z^k,
\end{equation}
and
\begin{equation}\label{defhypb}
_r\psi_s\!\left[\begin{matrix}a_1,a_2,\dots,a_r\\
b_1,b_2,\dots,b_s\end{matrix};q,z\right]:=
\sum _{k=-\infty} ^{\infty}\frac {(a_1,a_2,\dots,a_r;q)_k}
{(b_1,b_2,\dots,b_s;q)_k}\left((-1)^kq^{\binom k2}\right)^{s-r}z^k,
\end{equation}
for {\em basic hypergeometric $_r\phi_s$ series}, and {\em bilateral basic
hypergeometric $_r\psi_s$
series}, respectively. See \cite[p.~25 and p.~125]{grhyp} for the criteria
of when these series terminate, or, if not, when they converge. 
Note that, when we are considering  $_r\phi_s$ series
in this article, we usually have $s=r-1$, in which case the factor
$((-1)^kq^{\binom k2})^{1+s-r}$ in the series is just one.
A similar fact holds for the $_r\psi_r$ series when $s=r$.

To shorten some of our displays, we use Sears'~\cite{searstf} ``idem"
notation. The symbol ``$\operatorname{idem}(a;a_1,a_2,\dots,a_t)$"
after an expression stands for the sum of the $t$ expressions obtained
from the preceding expression by interchanging $a$ with each
$a_k$, $k=1,2,\dots,t$.

The theory of (classical) basic hypergeometric series consists of
several summation and transformation formulae involving $_{r+1}\phi_r$
or $_r\psi_r$ series. Some of the classical summation theorems require
that the parameters satisfy the additional condition of being
either balanced and/or very-well-poised.
An $_{r+1}\phi_r$ basic hypergeometric series is called
{\em balanced} if $b_1\cdots b_r=a_1\cdots a_{r+1}q$ and $z=q$.
An $_{r+1}\phi_r$ series is {\em well-poised} if
$a_1q=a_2b_1=\cdots=a_{r+1}b_r$. An $_{r+1}\phi_r$ basic
hypergeometric series is called {\em very-well-poised}
if it is well-poised and if $a_2=-a_3=q\sqrt{a_1}$.
Note that  the factor
\begin{equation}\label{vwp}
\frac {1-a_1q^{2k}}{1-a_1}
\end{equation}
appears in  a very-well-poised series.
The parameter $a_1$ is usually referred to as the
{\em special parameter} of such a series.
Similarly, a bilateral $_r\psi_r$ basic
hypergeometric series is well-poised if
$a_1b_1=a_2b_2\cdots=a_rb_r$ and very-well-poised if, in addition,
$a_1=-a_2=qb_1=-qb_2$.

In our computations in the following sections, we make heavily use of some
elementary identities involving $q$-shifted factorials which are listed
in Gasper and Rahman~\cite[Appendix~I]{grhyp}.

\section{Bailey's very-well-poised $_6\psi_6$ summation}\label{sec6}

To motivate some of our analysis in the later sections, we review here
our short and simple proof of Bailey's very-well-poised $_6\psi_6$ summation:
\begin{multline}\label{66gl}
{}_6\psi_6\!\left[\begin{matrix}q\sqrt{a},-q\sqrt{a},b,c,d,e\\
\sqrt{a},-\sqrt{a},aq/b,aq/c,aq/d,aq/e\end{matrix};q,
\frac{a^2q}{bcde}\right]\\
=\frac {(aq,aq/bc,aq/bd,aq/be,aq/cd,aq/ce,aq/de,q,q/a;q)_{\infty}}
{(aq/b,aq/c,aq/d,aq/e,q/b,q/c,q/d,q/e,a^2q/bcde;q)_{\infty}},
\end{multline}
provided the series either terminates, or $|a^2q/bcde|<1$,
for convergence.

The summation formula in \eqref{66gl} is one of the most powerful
identities for bilateral basic hypergeometric series. For some
applications to number theory, see Andrews~\cite[pp.~461--468]{andappl}.

To prove Bailey's $_6\psi_6$ summation, we start with a suitable
specialization of Rogers' $_6\phi_5$
summation:
\begin{equation}\label{65gl}
{}_6\phi_5\!\left[\begin{matrix}a,\,q\sqrt{a},-q\sqrt{a},b,c,d\\
\sqrt{a},-\sqrt{a},aq/b,aq/c,aq/d\end{matrix};q,
\frac{aq}{bcd}\right]
=\frac {(aq,aq/bc,aq/bd,aq/cd;q)_{\infty}}
{(aq/b,aq/c,aq/d,aq/bcd;q)_{\infty}},
\end{equation}
provided the series either terminates, or $|aq/bcd|<1$, for
convergence. Note that \eqref{65gl} is just the special case $e\mapsto a$
of \eqref{66gl}.

In \eqref{65gl}, we perform the simultaneous substitutions $a\mapsto c/a$,
$b\mapsto b/a$, $c\mapsto cq^{n}$ and $d\mapsto cq^{-n}/a$, and obtain
\begin{multline}\label{65gl1}
{}_6\phi_5\!\left[\begin{matrix}c/a,\,q\sqrt{c/a},-q\sqrt{c/a},b/a,cq^n,
cq^{-n}/a\\
\sqrt{c/a},-\sqrt{c/a},cq/b,q^{1-n}/a,q^{1+n}\end{matrix};q,
\frac{aq}{bc}\right]\\
=\frac {(cq/a,q^{1-n}/b,aq^{1+n}/b,q/c;q)_{\infty}}
{(cq/b,q^{1-n}/a,q^{1+n},aq/bc;q)_{\infty}},
\end{multline}
where $|aq/bc|<1$.

Using some elementary identities for $q$-shifted factorials
(see, e.g., Gasper and Rahman~\cite[Appendix~I]{grhyp})
we can rewrite equation~\eqref{65gl1} as
\begin{multline}\label{key1}
\frac {(cq/b,q/a,q,aq/bc;q)_{\infty}}
{(cq/a,q/b,aq/b,q/c;q)_{\infty}}
\sum_{k=0} ^{\infty}\frac {(1-cq^{2k}/a)} {(1-c/a)}
\frac {(c/a,b/a;q)_k(c;q)_{n+k}(a;q)_{n-k}}
{(q,cq/b;q)_k(q;q)_{n+k}(aq/c;q)_{n-k}}
\left(\frac {a} {b}\right)^k\\
=\frac {(b,c;q)_n}{(aq/b,aq/c;q)_n}\left(\frac {a} {b}\right)^n.
\end{multline}
In this identity, we multiply both sides by
\begin{equation*}
\frac{(1-aq^{2n})}{(1-a)}\frac{(d,e;q)_n}{(aq/d,aq/e;q)_n}
\left(\frac {aq} {cde}\right)^n
\end{equation*}
and sum over all integers $n$.

On the right side we obtain
\begin{equation*}
{}_6\psi_6\!\left[\begin{matrix}q\sqrt{a},-q\sqrt{a},b,c,d,e\\
\sqrt{a},-\sqrt{a},aq/b,aq/c,aq/d,aq/e\end{matrix};q,
\frac{a^2q}{bcde}\right].
\end{equation*}
On the left side we obtain
\begin{multline}\label{65gl2}
\frac {(cq/b,q/a,q,aq/bc;q)_{\infty}}
{(cq/a,q/b,aq/b,q/c;q)_{\infty}}
\sum_{n=-\infty} ^{\infty}\frac {(1-aq^{2n})} {(1-a)}
\frac {(d,e;q)_n} {(aq/d,aq/e;q)_n}
\left(\frac {aq} {cde}\right)^n\\
\times\sum_{k=0} ^{\infty}\frac {(1-cq^{2k}/a)} {(1-c/a)}
\frac {(c/a,b/a;q)_k(c;q)_{n+k}(a;q)_{n-k}}
{(q,cq/b;q)_k(q;q)_{n+k}(aq/c;q)_{n-k}}
\left(\frac {a} {b}\right)^k.
\end{multline}
Next, we interchange summations in \eqref{65gl2} and shift
the inner index $n\mapsto n-k$.
(Observe that the sum over $n$ is terminated by the term
$(q;q)_{n+k}^{-1}$ from below.)  We obtain, again using some elementary
identities for $q$-shifted factorials,
\begin{multline*}
\frac {(cq/b,q/a,q,aq/bc;q)_{\infty}}
{(cq/a,q/b,aq/b,q/c;q)_{\infty}}
\sum_{k=0} ^{\infty}\frac {(1-cq^{2k}/a)} {(1-c/a)}
\frac {(c/a,b/a;q)_k}{(q,cq/b;q)_k}\\
\times\frac{(1-aq^{-2k})}{(1-a)}
\frac{(a;q)_{-2k}(d,e;q)_{-k}}{(aq/c;q)_{-2k}(aq/d,aq/e;q)_{-k}}
\left(\frac {cde} {bq}\right)^k\\
\times\sum_{n=0} ^{\infty}\frac {(1-aq^{-2k+2n})} {(1-aq^{-2k})}
\frac {(aq^{-2k},c,dq^{-k},eq^{-k};q)_n}
{(q,aq^{1-2k}/c,aq^{1-k}/d,aq^{1-k}/e;q)_n}
\left(\frac {aq} {cde}\right)^n.
\end{multline*}
Now the inner sum, provided $|aq/cde|<1$,
can be evaluated by \eqref{65gl} and we obtain
\begin{multline*}
\frac {(cq/b,q/a,q,aq/bc;q)_{\infty}}
{(cq/a,q/b,aq/b,q/c;q)_{\infty}}
\sum_{k=0} ^{\infty}\frac {(1-cq^{2k}/a)} {(1-c/a)}
\frac {(c/a,b/a;q)_k(aq;q)_{-2k}}
{(q,cq/b;q)_k(aq/c;q)_{-2k}}\\
\times\frac{(d,e;q)_{-k}}{(aq/d,aq/e;q)_{-k}}
\left(\frac {cde} {bq}\right)^k
\frac{(aq^{1-2k},aq^{1-k}/cd,aq^{1-k}/ce,aq/de;q)_{\infty}}
{(aq^{1-2k}/c,aq^{1-k}/d,aq^{1-k}/e,aq/cde;q)_{\infty}},
\end{multline*}
which can be simplified to
\begin{multline*}
\frac {(cq/b,q/a,q,aq/bc,aq,aq/cd,aq/ce,aq/de;q)_{\infty}}
{(cq/a,q/b,aq/b,q/c,aq/c,aq/d,aq/e,aq/cde;q)_{\infty}}\\
\times\sum_{k=0} ^{\infty}\frac {(1-cq^{2k}/a)} {(1-c/a)}
\frac {(c/a,b/a,cd/a,ce/a;q)_k}
{(q,cq/b,q/d,q/e;q)_k}
\left(\frac {a^2q} {bcde}\right)^k.
\end{multline*}
To the last sum, provided $|a^2q/bcde|<1$, we can again apply
\eqref{65gl} and after some simplifications we finally obtain
the right side of \eqref{66gl}, as desired.

Our derivation of the $_6\psi_6$ summation~\eqref{66gl} is simple once
the nonterminating $_6\phi_5$ summation~\eqref{65gl} is given.
But the latter summation follows by an elementary computation from
F.~H.~Jackson's~\cite{jacksum} terminating
$_8\phi_7$ summation (cf.~\cite[Eq.~(2.6.2)]{grhyp})
\begin{multline}\label{87gl}
_8\phi_7\!\left[\begin{matrix}a,\,q\sqrt{a},-q\sqrt{a},b,c,d,
a^2q^{1+n}/bcd,q^{-n}\\
\sqrt{a},-\sqrt{a},aq/b,aq/c,aq/d,bcdq^{-n}/a,aq^{1+n}\end{matrix};q,
q\right]\\
=\frac {(aq,aq/bc,aq/bd,aq/cd;q)_n}
{(aq/b,aq/c,aq/d,aq/bcd;q)_n}
\end{multline}
as $n\to\infty$. Jackson's terminating $_8\phi_7$ summation itself can be
proved by various ways. An algorithmic approach uses the $q$-Zeilberger
algorithm, see Koornwinder~\cite{koornw}. For an inductive proof, see
Slater~\cite[Sec.~3.3.1]{slater}. For another elementary classical proof,
see Gasper and Rahman~\cite[Sec.~2.6]{grhyp}.

Concluding this section, we would like to add another thought,
kindly initiated by an anonymous referee of \cite{schlelsum}.
It is worth comparing our proof with Askey and Ismail's~\cite{askmail}
elegant (and now classical) proof of Bailey's $_6\psi_6$ summation.
Their proof uses a method in this context often referred to as
``Ismail's argument"
since Ismail~\cite{ismail} was apparently the first to apply
Liouville's standard analytic continuation argument in the context of
bilateral basic hypergeometric series.
Askey and Ismail use Rogers' $_6\phi_5$ summation once to evaluate
the $_6\psi_6$ series at an infinite sequence and then apply
analytic continuation.
Here, we evaluate the $_6\psi_6$ series on a domain, 
and, for the full theorem, we actually also need analytic continuation.
In fact, we need, in addition to $|a^2q/bcde|<1$ two other inequalities
on $a,b,c,d,e$, namely $|aq/bc|<1$ and $|aq/cde|<1$, in order to apply 
the $_6\phi_5$ summation theorem. In the end, these additional conditions
can be removed. In particular, both sides of identity \eqref{66gl}
are analytic in $1/c$ around the origin. So far, we have shown
the identity for $|1/c|<\min(|b/aq|,|de/aq|,|bde/a^2q|)$.
By analytic continuation, we extend the identity to be valid
for $|1/c|<|bde/a^2q|$, the radius of convergence of the series.

In the following sections, our objective is to find elementary
derivations of some of the classical {\em transformations} for bilateral
basic hypergeometric series. The method we will use is very similar to the
one used in this section with the difference that for the bilateral
transformations we also utilize bilateral series identities in our derivations.

\section{M.~Jackson's very-well-poised $_8\psi_8$ transformation}\label{sec8}

M.~Jackson's~\cite[Eq.~(2.2)]{mjack1} transformation
formula (cf.~\cite[Eq.~(5.6.2)]{grhyp}) of a very-well-poised
$_8\psi_8$ series into a sum of two $_8\phi_7$ series
can be stated as follows:
\begin{multline}\label{88gl}
{}_8\psi_8\!\left[\begin{matrix}q\sqrt{a},-q\sqrt{a},b,c,d,e,f,g\\
\sqrt{a},-\sqrt{a},aq/b,aq/c,aq/d,aq/e,aq/f,aq/g\end{matrix}\,;q,
\frac{a^3q^2}{bcdefg}\right]\\
=\frac{(q,aq,q/a,c,c/a,bq/d,bq/e,bq/f,bq/g,aq/bd,aq/be,aq/bf,
aq/bg;q)_{\infty}}
{(q/b,q/d,q/e,q/f,q/g,aq/b,aq/d,aq/e,aq/f,aq/g,c/b,bc/a,
b^2q/a;q)_{\infty}}\\\times
{}_8\phi_7\!\left[\begin{matrix}b^2/a,\,qb/\sqrt{a},-qb/\sqrt{a},
bc/a,bd/a,be/a,bf/a,bg/a\\
b/\sqrt{a},-b/\sqrt{a},bq/c,bq/d,bq/e,bq/f,bq/g\end{matrix}\,;q,
\frac{a^3q^2}{bcdefg}\right]\\
+\operatorname{idem}(b;c),
\end{multline}
where the series either terminate, or $|a^3q^2/bcdefg|<1$, for convergence.
(The standard symbol ``$\operatorname{idem}(b;c)$" is explained in
the introduction.)

M.~Jackson obtained this $_8\psi_8$ transformation formula by
specializing a general transformation of Sears~\cite{searstf}.

To derive the above $_8\psi_8$ transformation with our elementary method,
we start with a suitable specialization of Bailey's summation formula
for a nonterminating
very-well-poised $_8\phi_7$ series (cf.~\cite[Eq.~(2.11.7)]{grhyp}).
In our subsequent computations, we further make use of Bailey's
$_6\psi_6$ summation~\eqref{66gl}, and of Bailey's $_8\phi_7$ transformation
in \eqref{87tntgl}.

For convenience, we state Bailey's nonterminating
very-well-poised $_8\phi_7$ summation:
\begin{multline}\label{87ntgl}
{}_8\phi_7\!\left[\begin{matrix}a,\,q\sqrt{a},-q\sqrt{a},b,c,d,e,f\\
\sqrt{a},-\sqrt{a},aq/b,aq/c,aq/d,aq/e,aq/f\end{matrix}\,;q,q\right]\\
+\frac{(aq,c,d,e,f,b/a,bq/c,bq/d,bq/e,bq/f;q)_{\infty}}
{(a/b,aq/c,aq/d,aq/e,aq/f,bc/a,bd/a,be/a,bf/a,b^2q/a;q)_{\infty}}\\\times
{}_8\phi_7\!\left[\begin{matrix}b^2/a,\,qb/\sqrt{a},-qb/\sqrt{a},
b,bc/a,bd/a,be/a,bf/a\\
b/\sqrt{a},-b/\sqrt{a},bq/a,bq/c,bq/d,bq/e,bq/f\end{matrix}\,;q,q\right]\\
=\frac{(aq,b/a,aq/cd,aq/ce,aq/cf,aq/de,aq/df,aq/ef;q)_{\infty}}
{(aq/c,aq/d,aq/e,aq/f,bc/a,bd/a,be/a,bf/a;q)_{\infty}},
\end{multline}
where $a^2q=bcdef$. Let us briefly illustrate some of the depth
of this identity. Clearly, it contains the terminating
$_8\phi_7$ summation~\eqref{87gl} and nonterminating $_6\phi_5$
summation~\eqref{65gl} as special cases. A probably less known fact
is that the $b\to 1$ specialization of \eqref{87ntgl} yields an
important theta function identity (namely \cite[Eq.~(5.2)]{bail66}).
Of course, it is interesting to know how such a rich identity,
as \eqref{87ntgl} is, can be actually derived. Surprisingly, it can be
(not trivially) derived starting from a polynomial identity,
namely Bailey's~\cite{bail10} very-well-poised
$_{10}\phi_9$ transformation, see the exposition
in Gasper and Rahman~\cite[Secs.~2.10 and 2.11]{grhyp}.

To derive M.~Jackson's $_8\psi_8$ transformation, we perform
in \eqref{87ntgl} the simultaneous substitutions $a\mapsto abq/cde$,
$b\mapsto aq/de$, $c\mapsto aq/ce$, $d\mapsto bq^n$, and
$e\mapsto bq^{-n}/a$, and obtain
\begin{multline*}
{}_8\phi_7\!\left[\begin{matrix}abq/cde,\,q\sqrt{abq/cde},-q\sqrt{abq/cde},
aq/de,aq/ce,bq^n,bq^{-n}/a,aq/cd\\
\sqrt{abq/cde},-\sqrt{abq/cde},bq/c,bq/d,
aq^{2-n}/cde,a^2q^{2+n}/cde,bq/e\end{matrix}\,;q,q\right]\\
+\frac{(abq^2/cde,aq/ce,bq^n,bq^{-n}/a,aq/cd,c/b,cq/d,aq^{2-n}/bde,
a^2q^{2+n}/bde,cq/e;q)_{\infty}}
{(b/c,bq/d,aq^{2-n}/cde,a^2q^{2+n}/cde,bq/e,aq/be,cq^n,cq^{-n}/a,
aq/bd,acq^2/bde;q)_{\infty}}\\\times
{}_8\phi_7\!\left[\begin{matrix}acq/bde,\,q\sqrt{acq/bde},-q\sqrt{acq/bde},
aq/de,aq/be,cq^n,cq^{-n}/a,aq/bd\\
\sqrt{acq/bde},-\sqrt{acq/bde},cq/b,cq/d,aq^{2-n}/bde,a^2q^{2+n}/bde,
cq/e\end{matrix}\,;q,q\right]\\
=\frac{(abq^2/cde,c/b,q^{1-n}/d,aq^{1+n}/d,bc/a,a^2q^2/bcde,q^{1-n}/e,
aq^{1+n}/e;q)_{\infty}}
{(bq/d,aq^{2-n}/cde,a^2q^{2+n}/cde,bq/e,aq/be,cq^n,cq^{-n}/a,
aq/bd;q)_{\infty}}.
\end{multline*}
Using some elementary identities for $q$-shifted factorials
we can rewrite this as
\begin{multline}\label{87ntgl2}
\frac{(bq/d,aq^2/cde,a^2q^2/cde,bq/e,aq/be,c,c/a,aq/bd;q)_{\infty}}
{(abq^2/cde,c/b,q/d,aq/d,bc/a,a^2q^2/bcde,q/e,aq/e;q)_{\infty}}\\\times
\sum_{k=0} ^{\infty}\frac {(1-abq^{1+2k}/cde)} {(1-abq/cde)}
\frac {(abq/cde,aq/de,aq/ce,aq/cd;q)_k}
{(q,bq/c,bq/d,bq/e;q)_k}\\\times
\frac {(b;q)_{n+k}(cde/aq;q)_{n-k}}
{(a^2q^2/cde;q)_{n+k}(aq/b;q)_{n-k}}
\left(\frac {bcde} {a^2q}\right)^k\\
+\frac{(cq/d,aq^2/bde,a^2q^2/bde,cq/e,aq/ce,b,b/a,aq/cd;q)_{\infty}}
{(acq^2/bde,b/c,q/d,aq/d,bc/a,a^2q^2/bcde,q/e,aq/e;q)_{\infty}}\\\times
\sum_{k=0} ^{\infty}\frac {(1-acq^{1+2k}/bde)} {(1-acq/bde)}
\frac {(acq/bde,aq/de,aq/be,aq/bd;q)_k}
{(q,cq/b,cq/d,cq/e;q)_k}\\\times
\frac {(c;q)_{n+k}(bde/aq;q)_{n-k}}
{(a^2q^2/bde;q)_{n+k}(aq/c;q)_{n-k}}\left(\frac {bcde} {a^2q}\right)^k\\
=\frac {(b,c,d,e;q)_n}{(aq/b,aq/c,aq/d,aq/e;q)_n}.
\end{multline}

Observe that on the left side of this identity
the second term equals the first term where $b$ and $c$ are interchanged.
This observation helps us to reduce the amount of our subsequent
computations. 

In identity~\eqref{87ntgl2}, we multiply both sides by
\begin{equation*}
\frac{(1-aq^{2n})}{(1-a)}\frac{(f,g;q)_n}{(aq/f,aq/g;q)_n}
\left(\frac {a^3q^2} {bcdefg}\right)^n
\end{equation*}
and sum over all integers $n$.

On the right side we obtain
\begin{equation*}
{}_8\psi_8\!\left[\begin{matrix}q\sqrt{a},-q\sqrt{a},b,c,d,e,f,g\\
\sqrt{a},-\sqrt{a},aq/b,aq/c,aq/d,aq/e,aq/f,aq/g\end{matrix}\,;q,
\frac{a^3q^2}{bcdefg}\right].
\end{equation*}
On the left side we obtain
\begin{multline}\label{87ntgl3}
\frac{(bq/d,aq^2/cde,a^2q^2/cde,bq/e,aq/be,c,c/a,aq/bd;q)_{\infty}}
{(abq^2/cde,c/b,q/d,aq/d,bc/a,a^2q^2/bcde,q/e,aq/e;q)_{\infty}}\\\times
\sum_{n=-\infty} ^{\infty}\frac {(1-aq^{2n})} {(1-a)}
\frac{(f,g;q)_n}{(ag/f,aq/g;q)_n}\left(\frac{a^3q^2}{bcdefg}\right)^n\\\times
\sum_{k=0} ^{\infty}\frac {(1-abq^{1+2k}/cde)} {(1-abq/cde)}
\frac {(abq/cde,aq/de,aq/ce,aq/cd;q)_k}
{(q,bq/c,bq/d,bq/e;q)_k}\\\times
\frac {(b;q)_{n+k}(cde/aq;q)_{n-k}}
{(a^2q^2/cde;q)_{n+k}(aq/b;q)_{n-k}}
\left(\frac {bcde} {a^2q}\right)^k
+\operatorname{idem}(b;c).
\end{multline}
Next, we interchange summations in \eqref{87ntgl3} and obtain, again using
some elementary identities for $q$-shifted factorials,
\begin{multline*}
\frac{(bq/d,aq^2/cde,a^2q^2/cde,bq/e,aq/be,c,c/a,aq/bd;q)_{\infty}}
{(abq^2/cde,c/b,q/d,aq/d,bc/a,a^2q^2/bcde,q/e,aq/e;q)_{\infty}}\\\times
\sum_{k=0} ^{\infty}\frac {(1-abq^{1+2k}/cde)} {(1-abq/cde)}
\frac {(abq/cde,aq/de,aq/ce,aq/cd,b,b/a;q)_k}
{(q,bq/c,bq/d,bq/e,aq^2/cde,a^2q^2/cde;q)_k}q^k\\\times
\sum_{n=-\infty} ^{\infty}\frac {(1-aq^{2n})} {(1-a)}
\frac{(f,g,bq^k,cdeq^{-1-k}/a;q)_n}
{(ag/f,aq/g,aq^{1-k}/b,a^2q^{2+k}/cde;q)_n}
\left(\frac{a^3q^2}{bcdefg}\right)^n\\
+\operatorname{idem}(b;c).
\end{multline*}
Now the inner sums can be evaluated by the $_6\psi_6$
summation~\eqref{66gl} and we obtain
\begin{multline*}
\frac{(bq/d,aq^2/cde,a^2q^2/cde,bq/e,aq/be,c,c/a,aq/bd;q)_{\infty}}
{(abq^2/cde,c/b,q/d,aq/d,bc/a,a^2q^2/bcde,q/e,aq/e;q)_{\infty}}\\\times
\sum_{k=0} ^{\infty}\frac {(1-abq^{1+2k}/cde)} {(1-abq/cde)}
\frac {(abq/cde,aq/de,aq/ce,aq/cd,b,b/a;q)_k}
{(q,bq/c,bq/d,bq/e,aq^2/cde,a^2q^2/cde;q)_k}q^k\\\times
\frac{(aq,aq/fg,aq^{1-k}/bf,a^2q^{2+k}/cdef,aq^{1-k}/bg,a^2q^{2+k}/cdeg,
a^2q^2/bcde,q,q/a;q)_{\infty}}
{(aq/f,aq/g,aq^{1-k}/b,a^2q^{2+k}/cde,q/f,q/g,q^{1-k}/b,aq^{2+k}/cde,
a^3q^2/bcdefg;q)_{\infty}}\\
+\operatorname{idem}(b;c),
\end{multline*}
which can be simplified to
\begin{multline}\label{87ntgl4}
\frac{(bq/d,bq/e,aq/be,c,c/a,aq/bd;q)_{\infty}}
{(abq^2/cde,c/b,q/d,aq/d,bc/a,q/e,aq/e;q)_{\infty}}\\\times
\frac{(aq,aq/fg,aq/bf,a^2q^2/cdef,aq/bg,a^2q^2/cdeg,q,q/a;q)_{\infty}}
{(aq/f,aq/g,aq/b,q/f,q/g,q/b,a^3q^2/bcdefg;q)_{\infty}}\\\times
\sum_{k=0} ^{\infty}\frac {(1-abq^{1+2k}/cde)} {(1-abq/cde)}
\frac {(abq/cde,aq/de,aq/ce,aq/cd,bf/a,bg/a;q)_k}
{(q,bq/c,bq/d,bq/e,a^2q^2/cdef,a^2q^2/cdeg;q)_k}
\left(\frac{aq}{ef}\right)^k\\
+\operatorname{idem}(b;c).
\end{multline}

Now, to the $_8\phi_7$'s in~\eqref{87ntgl4} we apply
Bailey's~\cite[Eq.~(4.3)]{bail66} transformation formula
for a nonterminating $_8\phi_7$ series (see~\cite[Eq.~(2.10.1)]{grhyp}):
\begin{multline}\label{87tntgl}
{}_8\phi_7\!\left[\begin{matrix}a,\,q\sqrt{a},-q\sqrt{a},b,c,d,e,f\\
\sqrt{a},-\sqrt{a},aq/b,aq/c,aq/d,aq/e,aq/f\end{matrix}\,;q,
\frac{a^2q^2}{bcdef}\right]\\
=\frac{(aq,aq/ef,\lambda q/e,\lambda q/f;q)_{\infty}}
{(aq/e,aq/f,\lambda q,\lambda q/ef;q)_{\infty}}\\\times
{}_8\phi_7\!\left[\begin{matrix}\lambda,\,q\sqrt{\lambda},-q\sqrt{\lambda},
\lambda b/a,\lambda c/a,\lambda d/a,e,f\\
\sqrt{\lambda},-\sqrt{\lambda},aq/b,aq/c,aq/d,
\lambda q/e,\lambda q/f\end{matrix}\,;q,\frac{aq}{ef}\right],
\end{multline}
where $\lambda=a^2q/bcd$ and $\max(|aq/ef|, |\lambda q/ef|)<1$.
The purpose of our application of this transformation is to obtain
more symmetry.
Hence, by \eqref{87tntgl} the expression in \eqref{87ntgl4}
is transformed into
\begin{multline}\label{87ntgl5}
\frac{(bq/d,bq/e,aq/be,c,c/a,aq/bd,aq,aq/bf,aq/bg,q,q/a,
bq/f,bq/g;q)_{\infty}}
{(c/b,q/d,aq/d,bc/a,q/e,aq/e,aq/f,aq/g,aq/b,
q/f,q/g,q/b,b^2q/a;q)_{\infty}}\\\times
{}_8\phi_7\!\left[\begin{matrix}b^2/a,\,qb/\sqrt{\smash[b]a},
-qb/\sqrt{\smash[b]a},bc/a,bd/a,be/a,bf/a,bg/a\\
b/\sqrt{\smash[b]a},-b/\sqrt{\smash[b]a},bq/c,bq/d,bq/e,bq/f,
bq/g\end{matrix}\,;q,
\frac{a^3q^2}{bcdefg}\right]\\
+\operatorname{idem}(b;c),
\end{multline}
which already completes our derivation of M.~Jackson's very-well-poised
$_8\psi_8$ transformation formula.

\section{A very-well-poised $_{10}\psi_{10}$ transformation}\label{sec10}

In this section we derive with our method the following
very-well-poised $_{10}\psi_{10}$ transformation 
(cf.~\cite[Eq.~(5.6.3)]{grhyp}):
\begin{multline}\label{1010gl}
{}_{10}\psi_{10}\!\left[\begin{matrix}q\sqrt{a},-q\sqrt{a},b,c,d,e,f,g,h,y\\
\sqrt{a},-\sqrt{a},\frac{aq}b,\frac{aq}c,\frac{aq}d,\frac{aq}e,\frac{aq}f,
\frac{aq}g,\frac{aq}h,\frac{aq}y\end{matrix}\,;q,
\frac{a^4q^3}{bcdefghy}\right]\\
=\frac{(q,aq,\frac qa,c,\frac ca,d,\frac da,\frac{bq}e,\frac{bq}f,\frac{bq}g,
\frac{bq}h,\frac{bq}y,\frac{aq}{be},\frac{aq}{bf},\frac{aq}{bg},
\frac{aq}{bh},\frac{aq}{by};q)_{\infty}}
{(\frac qb,\frac qe,\frac qf,\frac qg,\frac qh,\frac qy,
\frac{aq}b,\frac{aq}e,\frac{aq}f,\frac{aq}g,\frac{aq}h,\frac{aq}y,
\frac{bc}a,\frac{bd}a,\frac cb,\frac db,\frac{b^2q}a;q)_{\infty}}\\\times
{}_{10}\phi_9\!\left[\begin{matrix}\frac{b^2}a,\,
\frac{qb}{\sqrt{a}},-\frac {qb}{\sqrt{a}},
\frac{bc}a,\frac{bd}a,\frac{be}a,\frac{bf}a,\frac{bg}a,
\frac{bh}a,\frac{by}a\\
\frac b{\sqrt{a}},-\frac b{\sqrt{a}},
\frac{bq}c,\frac{bq}d,\frac{bq}e,\frac{bq}f,\frac{bq}g,
\frac{bq}h,\frac{bq}y\end{matrix}\,;q,
\frac{a^4q^3}{bcdefghy}\right]\\
+\operatorname{idem}(b;c,d),
\end{multline}
where the series either terminate, or $|a^4q^3/bcdefghy|<1$, for convergence.
(The symbol ``$\operatorname{idem}(b;c,d)$" is explained in the introduction.)

Gasper and Rahman~\cite[Sec.~5.6]{grhyp} derive this $_{10}\psi_{10}$
transformation formula by specializing a general
transformation of Slater~\cite{slatertf}.

To derive this transformation formula with our method, we make use of
a special case of Bailey's nonterminating $_8\phi_7$ summation
\eqref{87ntgl}, namely the key identity \eqref{87ntgl2},
and of M.~Jackson's $_8\psi_8$ summation \eqref{88gl} which we just derived
in the previous section. In the course of our derivation, we also make
use of Jackson's terminating $_8\phi_7$ summation \eqref{87gl}, and of
Bailey's nonterminating $_8\phi_7$ summation \eqref{87ntgl}
in its general form.

In identity~\eqref{87ntgl2}, we first replace $d$ by $f$. Then we
multiply both sides by
\begin{equation*}
\frac{(1-aq^{2n})}{(1-a)}\frac{(d,g,h,y;q)_n}{(aq/d,aq/g,aq/h,aq/y;q)_n}
\left(\frac {a^4q^3} {bcdefghy}\right)^n
\end{equation*}
and sum over all integers $n$.

On the right side we obtain
\begin{equation*}
{}_{10}\psi_{10}\!\left[\begin{matrix}q\sqrt{a},-q\sqrt{a},b,c,d,e,f,g,h,y\\
\sqrt{a},-\sqrt{a},aq/b,aq/c,aq/d,aq/e,aq/f,aq/g,aq/h,aq/y\end{matrix}\,;q,
\frac{a^4q^3}{bcdefghy}\right].
\end{equation*}
On the left side we obtain
\begin{multline}\label{1010ntgl1}
\frac{(bq/f,aq^2/cef,a^2q^2/cef,bq/e,aq/be,c,c/a,aq/bf;q)_{\infty}}
{(abq^2/cef,c/b,q/f,aq/f,bc/a,a^2q^2/bcef,q/e,aq/e;q)_{\infty}}\\\times
\sum_{n=-\infty} ^{\infty}\frac {(1-aq^{2n})} {(1-a)}
\frac{(d,g,h,y;q)_n}{(aq/d,aq/g,ag/h,aq/y;q)_n}
\left(\frac{a^4q^3}{bcdefghy}\right)^n\\\times
\sum_{k=0} ^{\infty}\frac {(1-abq^{1+2k}/cef)} {(1-abq/cef)}
\frac {(abq/cef,aq/ef,aq/ce,aq/cf;q)_k}
{(q,bq/c,bq/f,bq/e;q)_k}\\\times
\frac {(b;q)_{n+k}(cef/aq;q)_{n-k}}
{(a^2q^2/cef;q)_{n+k}(aq/b;q)_{n-k}}
\left(\frac {bcef} {a^2q}\right)^k
+\operatorname{idem}(b;c).
\end{multline}
Next, we interchange summations in \eqref{1010ntgl1} and obtain
\begin{multline*}
\frac{(bq/f,aq^2/cef,a^2q^2/cef,bq/e,aq/be,c,c/a,aq/bf;q)_{\infty}}
{(abq^2/cef,c/b,q/f,aq/f,bc/a,a^2q^2/bcef,q/e,aq/e;q)_{\infty}}\\\times
\sum_{k=0} ^{\infty}\frac {(1-abq^{1+2k}/cef)} {(1-abq/cef)}
\frac {(abq/cef,aq/ef,aq/ce,aq/cf,b,b/a;q)_k}
{(q,bq/c,bq/f,bq/e,aq^2/cef,a^2q^2/cef;q)_k}q^k\\\times
\sum_{n=-\infty} ^{\infty}\frac {(1-aq^{2n})} {(1-a)}
\frac{(bq^k,d,g,h,y,cefq^{-1-k}/a;q)_n}
{(aq^{1-k}/b,aq/d,ag/q,aq/h,aq/y,a^2q^{2+k}/cef;q)_n}
\left(\frac{a^4q^3}{bcdefghy}\right)^n\\
+\operatorname{idem}(b;c).
\end{multline*}
Now to the inner sums we apply M.~Jackson's $_8\psi_8$ transformation
\eqref{88gl} and we obtain
\begin{multline*}
\frac{(bq/f,aq^2/cef,a^2q^2/cef,bq/e,aq/be,c,c/a,aq/bf;q)_{\infty}}
{(abq^2/cef,c/b,q/f,aq/f,bc/a,a^2q^2/bcef,q/e,aq/e;q)_{\infty}}\\\times
\sum_{k=0} ^{\infty}\frac {(1-abq^{1+2k}/cef)} {(1-abq/cef)}
\frac {(abq/cef,aq/ef,aq/ce,aq/cf,b,b/a;q)_k}
{(q,bq/c,bq/f,bq/e,aq^2/cef,a^2q^2/cef;q)_k}q^k\\\times
\Bigg[
\frac{(q,aq,q/a,d,d/a,bq^{1+k}/g,bq^{1+k}/h,bq^{1+k}/y,
abq^{2+2k}/cef,aq^{1-k}/bg;q)_{\infty}}
{(q^{1-k}/b,q/g,q/h,q/y,aq^{2+k}/cef,aq^{1-k}/b,aq/g,
aq/h,aq/y,a^2q^{2+k}/cef;q)_{\infty}}\\\times
\frac{(aq^{1-k}/bh,aq^{1-k}/by,a^2q/bcef;q)_{\infty}}
{(dq^{-k}/b,bdq^k/a,b^2q^{1+2k}/a;q)_{\infty}}
\sum_{n=0}^{\infty}\frac{(1-b^2q^{2k+2n}/a)}{(1-b^2q^{2k}/a)}
\frac{(b^2q^{2k}/a,bdq^k/a;q)_n}
{(q,bq^{1+k}/d;q)_n}\\\times
\frac{(bgq^k/a,bhq^k/a,byq^k/a,bcef/a^2q;q)_n}
{(bq^{1+k}/g,bq^{1+k}/h,bq^{1+k}/y,abq^{2+2k}/cef;q)_n}
\left(\frac{a^4q^3}{bcdefghy}\right)^n\\
+\frac{(q,aq,q/a,bq^k,bq^k/a,dq/g,dq/h,dq/y,adq^{2+2k}/cef,aq/dg;q)_{\infty}}
{(q/d,q/g,q/h,q/y,aq^{2+k}/cef,aq/d,aq/g,aq/h,aq/y,
a^2q^{2+k}/cef;q)_{\infty}}\\\times
\frac{(aq/dh,aq/dy,a^2q^{2+k}/cdef;q)_{\infty}}
{(bq^k/d,bdq^k/a,d^2q/a;q)_{\infty}}
\sum_{n=0}^{\infty}\frac{(1-d^2q^{2n}/a)}{(1-d^2/a)}
\frac{(d^2/a,dbq^k/a;q)_n}
{(q,dq^{1-k}/b;q)_n}\\\times
\frac{(dg/a,dh/a,dy/a,cdefq^{-1-k}/a^2;q)_n}
{(dq/g,dq/h,dq/y,adq^{2+k}/cef;q)_n}
\left(\frac{a^4q^3}{bcdefghy}\right)^n
\Bigg]
+\operatorname{idem}(b;c),
\end{multline*}
which can be (slightly) simplified to
\begin{multline}\label{1010ntgl2}
\Bigg[\frac{(bq/f,bq/e,aq/be,c,c/a,aq/bf,q,aq,q/a,d,d/a,bq/g,bq/h,
bq/y;q)_{\infty}}
{(c/b,q/f,aq/f,bc/a,q/e,aq/e,q/b,q/g,q/h,q/y,aq/b,aq/g,
aq/h,aq/y;q)_{\infty}}\\\times
\frac{(aq/bg,aq/bh,aq/by;q)_{\infty}}{(d/b,bd/a,b^2q/a;q)_{\infty}}
\sum_{k=0} ^{\infty}\frac {(1-abq^{1+2k}/cef)} {(1-abq/cef)}
\frac {(abq/cef,aq/ef,aq/ce;q)_k}
{(q,bq/c,bq/f;q)_k}\\\times
\frac{(aq/cf,bd/a,bg/a,bh/a,by/a;q)_k(b^2q/a;q)_{2k}}
{(bq/e,bq/d,bq/g,bq/h,bq/y;q)_k(ab^2q^2/cef;q)_{2k}}
\left(\frac{a^2q^2}{dghy}\right)^k\\\times
\sum_{n=0}^{\infty}\frac{(1-b^2q^{2k+2n}/a)}{(1-b^2q^{2k}/a)}
\frac{(b^2q^{2k}/a,bdq^k/a,bgq^k/a,bhq^k/a;q)_n}
{(q,bq^{1+k}/d,bq^{1+k}/g,bq^{1+k}/h;q)_n}\\\times
\frac{(byq^k/a,bcef/a^2q;q)_n}
{(bq^{1+k}/y,abq^{2+2k}/cef;q)_n}
\left(\frac{a^4q^3}{bcdefghy}\right)^n\\
+\frac{(bq/f,bq/e,aq/be,c,c/a,aq/bf,q,aq,q/a,b,b/a,dq/g;q)_{\infty}}
{(abq^2/cef,c/b,q/f,aq/f,bc/a,a^2q^2/bcef,q/e,aq/e,q/d,q/g,q/h,
q/y;q)_{\infty}}\\\times
\frac{(dq/h,dq/y,adq^2/cef,aq/dg,aq/dh,aq/dy,a^2q^2/cdef;q)_{\infty}}
{(aq/d,aq/g,aq/h,aq/y,b/d,bd/a,d^2q/a;q)_{\infty}}\\\times
\sum_{k=0} ^{\infty}\frac {(1-abq^{1+2k}/cef)} {(1-abq/cef)}
\frac {(abq/cef,aq/ef,aq/ce,aq/cf,b/d,bd/a;q)_k}
{(q,bq/c,bq/f,bq/e,adq^2/cef,a^2q^2/cdef;q)_k}q^k\\\times
\sum_{n=0}^{\infty}\frac{(1-d^2q^{2n}/a)}{(1-d^2/a)}
\frac{(d^2/a,dbq^k/a,dg/a,dh/a;q)_n}
{(q,dq^{1-k}/b,dg/a,dh/a;q)_n}\\\times
\frac{(dy/a,cdefq^{-1-k}/a^2;q)_n}
{(dq/y,adq^{2+k}/cef;q)_n}
\left(\frac{a^4q^3}{bcdefghy}\right)^n
\Bigg]
+\operatorname{idem}(b;c).
\end{multline}

We have in \eqref{1010ntgl2} a sum of four double sums.
Accordingly, for more clarity, let us write the whole expression
in \eqref{1010ntgl2} as
\begin{equation}\label{t1234}
T_1(b,c)+T_2(b,c)+U_1(b,c)+U_2(b,c),
\end{equation}
where, by definition of ``idem", $U_1(b,c)=T_1(c,b)$ and $U_2(b,c)=T_2(c,b)$.
Below, we will selectively perform manipulations with the terms $T_1$,
$T_2$, $U_1$, and $U_2$. 

To evaluate $T_1(b,c)$ (and hence also $U_1(b,c)$),
we first shift the index $n$ of the inner sum in $T_1(b,c)$ (or,
equivalently, in the first term in \eqref{1010ntgl2}) by $-k$
and then interchange the double sum.
Symbolically, we apply
\begin{equation}\label{intchs}
\sum_{k=0}^{\infty}\sum_{n=0}^{\infty}f(n,k)=\sum_{n=0}^{\infty}\sum_{k=0}^n
f(n-k,k).
\end{equation}
Thus, we obtain (using some elementary identities for $q$-shifted factorials)
\begin{multline*}
T_1(b,c)=
\frac{(bq/f,bq/e,aq/be,c,c/a,aq/bf,q,aq,q/a,d,d/a;q)_{\infty}}
{(c/b,q/f,aq/f,bc/a,q/e,aq/e,q/b,q/g,q/h,q/y,aq/b;q)_{\infty}}\\\times
\frac{(bq/g,bq/h,bq/y,aq/bg,aq/bh,aq/by;q)_{\infty}}
{(aq/g,aq/h,aq/y,d/b,bd/a,b^2q/a;q)_{\infty}}\\\times
\sum_{n=0}^{\infty}\frac{(1-b^2q^{2n}/a)}{(1-b^2/a)}
\frac{(b^2/a,bd/a,bg/a,bh/a,by/a,bcef/a^2q;q)_n}
{(q,bq/d,bq/g,bq/h,bq/y,abq^2/cef;q)_n}
\left(\frac{a^4q^3}{bcdefghy}\right)^n\\\times
\sum_{k=0} ^n\frac {(1-abq^{1+2k}/cef)} {(1-abq/cef)}
\frac {(abq/cef,aq/ef,aq/ce,aq/cf,b^2q^n/a,q^{-n};q)_k}
{(q,bq/c,bq/f,bq/e,a^2q^{2-n}/bcef,abq^{2+n}/cef;q)_k}q^k.
\end{multline*}
Now the inner sum can be evaluated by Jackson's terminating
$_8\phi_7$ summation \eqref{87gl}, which simplifies the last
expression, $T_1(b,c)$, to 
\begin{multline}\label{1010ntgl3}
\frac{(bq/f,bq/e,aq/be,c,c/a,aq/bf,q,aq,q/a,d,d/a,bq/g,bq/h,
bq/y;q)_{\infty}}
{(c/b,q/f,aq/f,bc/a,q/e,aq/e,q/b,q/g,q/h,q/y,aq/b,aq/g,
aq/h,aq/y;q)_{\infty}}\\\times
\frac{(aq/bg,aq/bh,aq/by;q)_{\infty}}{(d/b,bd/a,b^2q/a;q)_{\infty}}
\sum_{n=0}^{\infty}\frac{(1-b^2q^{2n}/a)}{(1-b^2/a)}
\frac{(b^2/a,bc/a,bd/a,be/a;q)_n}
{(q,bq/c,bq/d,bq/e;q)_n}\\\times
\frac{(bf/a,bg/a,bh/a,by/a,;q)_n}
{(bq/f,bq/g,bq/h,bq/y;q)_n}
\left(\frac{a^4q^3}{bcdefghy}\right)^n.
\end{multline}

Next, we consider $T_2(b,c)$ and $U_2(b,c)$.
By interchanging the double sums in $T_2(b,c)$ and in $U_2(b,c)$ we obtain
\begin{multline}\label{1010ntgl4}
T_2(b,c)+U_2(b,c)=
\frac{(bq/f,bq/e,aq/be,c,c/a,aq/bf,q,aq,q/a;q)_{\infty}}
{(abq^2/cef,c/b,q/f,aq/f,bc/a,a^2q^2/bcef,q/e,aq/e,q/d;q)_{\infty}}\\\times
\frac{(b,b/a,dq/g,dq/h,dq/y,adq^2/cef,aq/dg,aq/dh,aq/dy,
a^2q^2/cdef;q)_{\infty}}
{(q/g,q/h,q/y,aq/d,aq/g,aq/h,aq/y,b/d,bd/a,d^2q/a;q)_{\infty}}\\\times
\sum_{n=0}^{\infty}\frac{(1-d^2q^{2n}/a)}{(1-d^2/a)}
\frac{(d^2/a,db/a,dg/a,dh/a,dy/a,cdef/a^2q;q)_n}
{(q,dq/b,dg/a,dh/a,dq/y,adq^2/cef;q)_n}
\left(\frac{a^4q^3}{bcdefghy}\right)^n\\\times
\sum_{k=0} ^{\infty}\frac {(1-abq^{1+2k}/cef)} {(1-abq/cef)}
\frac {(abq/cef,aq/ef,aq/ce,aq/cf,bq^{-n}/d,bdq^n/a;q)_k}
{(q,bq/c,bq/f,bq/e,adq^{2+n}/cef,a^2q^{2-n}/cdef;q)_k}q^k\\
+\operatorname{idem}(b;c).
\end{multline}
Now, to the inner sum of the {\em first} double sum (but not of the
second!) in \eqref{1010ntgl4} we apply
Bailey's nonterminating $_8\phi_7$ summation \eqref{87ntgl}, i.e.,
specifically we apply
\begin{multline}\label{87ntgl10}
\sum_{k=0}^{\infty}\frac{(1-abq^{1+2k}/cef)}{(1-abq/cef)}
\frac{(abq/cef,aq/ef,aq/ce,aq/cf,bq^{-n}/d,bdq^n/a;q)_k}
{(q,bq/c,bq/f,bq/e,adq^{2+n}/cef,a^2q^{2-n}/cdef;q)_k}q^k\\
=\frac{(abq^2/cef,c/b,bc/a,dq^{1+n}/f,aq^{1-n}/df,dq^{1+n}/e,aq^{1-n}/de,
a^2q^2/bcef;q)_{\infty}}
{(bq/f,bq/e,adq^{2+n}/cef,a^2q^{2-n}/cdef,aq/be,aq/bf,cq^{-n}/d,
cdq^n/a;q)_{\infty}}\\
-\frac{(abq^2/cef,aq/ce,aq/cf,bq^{-n}/d,bdq^n/a,c/b,cq/f,cq/e,
adq^{2+n}/bef;q)_{\infty}}
{(b/c,bq/f,bq/e,adq^{2+n}/cef,a^2q^{2-n}/cdef,aq/be,aq/bf,cq^{-n}/d,
cdq^n/a;q)_{\infty}}\\\times
\frac{(a^2q^{2-n}/bdef;q)_{\infty}}{(acq^2/bef;q)_{\infty}}
\sum_{k=0}^{\infty}\frac{(1-acq^{1+2k}/bef)}{(1-acq/bef)}
\frac{(acq/bef,aq/ef,aq/be;q)_k}
{(q,cq/b,cq/f;q)_k}\\\times
\frac{(aq/bf,cq^{-n}/d,cdq^n/a;q)_k}
{(cq/e,adq^{2+n}/bef,a^2q^{2-n}/bdef;q)_k}q^k.
\end{multline}
The result of the application of this (two term) summation is that
the first term in \eqref{1010ntgl4}, $T_2(b,c)$, is split into two parts,
one single sum $T_2^{\prime}(b,c)$ and one double sum $T_2^{\prime\prime}(b,c)$.
Formally, we have
$$
T_2(b,c)+U_2(b,c)=[T_2^{\prime}(b,c)+T_2^{\prime\prime}(b,c)]+U_2(b,c).
$$
But $T_2^{\prime\prime}(b,c)$ is precisely $-U_2(b,c)$
(as can be readily checked), so two terms cancel, thus
$$
T_2(b,c)+U_2(b,c)=T_2^{\prime}(b,c).
$$
Now, the last expression, $T_2^{\prime}(b,c)$, can be simplified to
\begin{multline}\label{1010ntgl5}
\frac{(c,c/a,q,aq,q/a,b,b/a,dq/g,dq/h,dq/y,aq/dg,aq/dh,aq/dy,dq/f;q)_{\infty}}
{(q/f,aq/f,q/e,aq/e,q/d,q/g,q/h,q/y,aq/d,aq/g,aq/h,aq/y,
b/d,bd/a;q)_{\infty}}\\\times
\frac{(aq/df,dq/e,aq/de;q)_{\infty}}
{(d^2q/a,c/d,cd/a;q)_{\infty}}
\sum_{n=0}^{\infty}\frac{(1-d^2q^{2n}/a)}{(1-d^2/a)}
\frac{(d^2/a,db/a,dc/a,de/a;q)_n}{(q,dq/b,dq/c,dq/e;q)_n}\\\times
\frac{(df/a,dg/a,dh/a,dy/a;q)_n}
{(dq/f,dg/a,dh/a,dq/y;q)_n}
\left(\frac{a^4q^3}{bcdefghy}\right)^n.
\end{multline}
It is easy to see that \eqref{1010ntgl5} equals
\eqref{1010ntgl3} where $b$ and $d$ are interchanged.
Collecting all terms, according to \eqref{t1234}, completes
our derivation of the very-well-poised $_{10}\psi_{10}$
transformation in \eqref{1010gl}.

\section{Slater's transformations for bilateral well-poised series}
\label{sec2r}

In view of the success of our (more or less) elaborate but elementary
derivation of the very-well-poised $_{10}\psi_{10}$ transformation
in the preceding section, we are encouraged to go for more.
In fact, the same machinery applies, together with induction, to
prove the following general very-well-poised $_{2r}\psi_{2r}$
transformation formula due to Slater~\cite{slatertf}.
For $r\ge 3$,
\begin{multline}\label{2r2rgl}
{}_{2r}\psi_{2r}\!\left[\begin{matrix}q\sqrt{a},-q\sqrt{a},b_3,b_4,
\dots,b_{2r}\\
\sqrt{a},-\sqrt{a},\frac{a q}{b_3},\frac{a q}{b_4},
\dots,\frac{a q}{b_{2r}}\end{matrix}\,;q,
\frac{a^{r-1}q^{r-2}}{b_3\dots b_{2r}}\right]\\
=\frac{(q,a q,\frac q{a},b_4,\dots,b_r,\frac{b_4}{a},\dots,\frac{b_r}{a},
\frac{b_3q}{b_{r+1}},\dots,\frac{b_3q}{b_{2r}},
\frac{a q}{b_3b_{r+1}},\dots,\frac{a q}{b_3b_{2r}};q)_{\infty}}
{(\frac q{b_3},\frac q{b_{r+1}},\dots,\frac q{b_{2r}},
\frac{a q}{b_3},\frac{a q}{b_{r+1}},\dots,\frac{a q}{b_{2r}},
\frac{b_4}{b_3},\dots,\frac{b_r}{b_3},\frac{b_3b_4}{a},\dots,
\frac{b_3b_r}{a},\frac{b_3^2q}{a};q)_{\infty}}\\\times
{}_{2r}\phi_{2r-1}\!\left[\begin{matrix}\frac{b_3^2}{a},\,
\frac{qb_3}{\sqrt{a}},-\frac {qb_3}{\sqrt{a}},
\frac{b_3b_4}{a},\frac{b_3b_5}{a},\dots,\frac{b_3b_{2r}}{a}\\
\frac{b_3}{\sqrt{a}},-\frac{b_3}{\sqrt{a}},
\frac{b_3q}{b_4},\frac{b_3q}{b_5},\dots,\frac{b_3q}{b_{2r}}\end{matrix}\,;q,
\frac{a^{r-1}q^{r-2}}{b_3\dots b_{2r}}\right]\\
+\operatorname{idem}(b_3;b_4,\dots,b_r),
\end{multline}
where the series either terminate, or $|a^{r-1}q^{r-2}/b_3\dots b_{2r}|<1$,
for convergence.
(The symbol ``$\operatorname{idem}(b_3;b_4,\dots,b_r)$" is explained in
the introduction.)

Slater~\cite{slatertf} first obtained this transformation formula
(or rather the more general one in \eqref{sl2r2rgl}) 
by extending Sears'~\cite{searstf} general
transformation for (unilateral) basic hypergeometric series.
Not much later, she~\cite{slatertf2} gave a direct proof using
a basic Barnes-type contour integral.

Here, we provide an inductive proof of \eqref{2r2rgl}
following closely the analysis of our derivation of the
very-well-poised $_{10}\psi_{10}$ transformation in Section~\ref{sec10}. 

The $r=3$ case of \eqref{2r2rgl} is readily
verified using the $_6\psi_6$ and $_6\phi_5$ summations in
\eqref{66gl} and in \eqref{65gl}, respectively.
We have also shown the $r=4$ and $r=5$ cases in Sections~\ref{sec8}
and \ref{sec10}. So we can assume that $r>5$.

Now suppose that the formula in \eqref{2r2rgl}
is already shown for any integer $t$ where $3\le t<r$.
To show the transformation for $t=r$,
we first perform in the key identity~\eqref{87ntgl2} the substitutions
$b\mapsto b_3$, $c\mapsto b_4$, $d\mapsto b_{r+2}$,
and $e\mapsto b_{r+1}$.
In the resulting equation, we multiply both sides by
\begin{equation*}
\frac{(1-a q^{2n})}{(1-a)}\frac{(b_5,\dots,b_r,b_{r+3},\dots,b_{2r};q)_n}
{(a q/b_5,\dots,a q/b_r,a q/b_{r+3},\dots,a q/b_{2r};q)_n}
\left(\frac{a^{r-1}q^{r-2}}{b_3\dots b_{2r}}\right)^n
\end{equation*}
and sum over all integers $n$.

On the right side we obtain
\begin{equation*}
{}_{2r}\psi_{2r}\!\left[\begin{matrix}q\sqrt{a},-q\sqrt{a},b_3,b_4,\dots,
b_{2r}\\
\sqrt{a},-\sqrt{a},a q/b_3,a q/b_4,\dots,a q/b_{2r}\end{matrix}\,;q,
\frac{a^{r-1}q^{r-2}}{b_3\dots b_{2r}}\right].
\end{equation*}
On the left side we obtain
\begin{multline}\label{2r2rgl1}
\frac{(b_3q/b_{r+2},a q^2/b_4b_{r+1}b_{r+2},a^2q^2/b_4b_{r+1}b_{r+2},
b_3q/b_{r+1};q)_{\infty}}
{(a b_3q^2/b_4b_{r+1}b_{r+2},b_4/b_3,q/b_{r+2},
a q/b_{r+2};q)_{\infty}}\\\times
\frac{(a q/b_3b_{r+1},b_4,b_4/a,a q/b_3b_{r+2};q)_{\infty}}
{(b_3b_4/a,a^2q^2/b_3b_4b_{r+1}b_{r+2},q/b_{r+1},
a q/b_{r+1};q)_{\infty}}\\\times
\sum_{n=-\infty} ^{\infty}\frac {(1-a q^{2n})} {(1-a)}
\frac{(b_5,\dots,b_r,b_{r+3},\dots,b_{2r};q)_n}
{(a q/b_5,\dots,a q/b_r,a q/b_{r+3},\dots,a q/b_{2r};q)_n}
\left(\frac{a^{r-1}q^{r-2}}{b_3\dots b_{2r}}\right)^n\\\times
\sum_{k=0} ^{\infty}\frac {(1-a b_3q^{1+2k}/b_4b_{r+1}b_{r+2})}
{(1-a b_3q/b_4b_{r+1}b_{r+2})}
\frac {(a b_3q/b_4b_{r+1}b_{r+2},a q/b_{r+1}b_{r+2},a q/b_4b_{r+1};q)_k}
{(q,b_3q/b_4,b_3q/b_{r+2};q)_k}\\\times
\frac {(a q/b_4b_{r+2};q)_k
(b_3;q)_{n+k}(b_4b_{r+1}b_{r+2}/a q;q)_{n-k}}
{(b_3q/b_{r+1};q)_k(a^2q^2/b_4b_{r+1}b_{r+2};q)_{n+k}
(a q/b_3;q)_{n-k}}
\left(\frac {b_3b_4b_{r+1}b_{r+2}} {a^2q}\right)^k\\
+\operatorname{idem}(b_3;b_4).
\end{multline}
Next, we interchange summations in \eqref{2r2rgl1} and obtain
\begin{multline}\label{2r2rgl2}
\frac{(b_3q/b_{r+2},a q^2/b_4b_{r+1}b_{r+2},a^2q^2/b_4b_{r+1}b_{r+2},
b_3q/b_{r+1};q)_{\infty}}
{(a b_3q^2/b_4b_{r+1}b_{r+2},b_4/b_3,q/b_{r+2},
a q/b_{r+2};q)_{\infty}}\\\times
\frac{(a q/b_3b_{r+1},b_4,b_4/a,a q/b_3b_{r+2};q)_{\infty}}
{(b_3b_4/a,a^2q^2/b_3b_4b_{r+1}b_{r+2},q/b_{r+1},
a q/b_{r+1};q)_{\infty}}\\\times
\sum_{k=0} ^{\infty}\frac {(1-a b_3q^{1+2k}/b_4b_{r+1}b_{r+2})}
{(1-a b_3q/b_4b_{r+1}b_{r+2})}
\frac {(a b_3q/b_4b_{r+1}b_{r+2},a q/b_{r+1}b_{r+2},a q/b_4b_{r+1};q)_k}
{(q,b_3q/b_4,b_3q/b_{r+2};q)_k}\\\times
\frac {(a q/b_4b_{r+2},b_3,b_3/a;q)_k}
{(b_3q/b_{r+1},a q^2/b_4b_{r+1}b_{r+2},a^2 q^2/b_4b_{r+1}b_{r+2};q)_k}q^k\\
\sum_{n=-\infty} ^{\infty}\frac {(1-a q^{2n})} {(1-a)}
\frac{(b_3q^k,b_5,\dots,b_r,b_{r+3},\dots,b_{2r};q)_n}
{(a q^{1-k}/b_3,a q/b_5,\dots,a q/b_r,a q/b_{r+3},
\dots,a q/b_{2r};q)_n}\\\times
\frac{(b_4b_{r+1}b_{r+2}q^{-1-k}/a;q)_n}
{(a^2q^{2+k}/b_4b_{r+1}b_{r+2};q)_n}
\left(\frac{a^{r-1}q^{r-2}}{b_3\dots b_{2r}}\right)^n
+\operatorname{idem}(b_3;b_4).
\end{multline}
Now, to the inner sums we apply the inductive hypothesis
(i.e., the $r\mapsto r-1$ case of \eqref{2r2rgl}), and we obtain
\begin{multline*}
\frac{(b_3q/b_{r+2},a q^2/b_4b_{r+1}b_{r+2},a^2q^2/b_4b_{r+1}b_{r+2},
b_3q/b_{r+1};q)_{\infty}}
{(a b_3q^2/b_4b_{r+1}b_{r+2},b_4/b_3,q/b_{r+2},
a q/b_{r+2};q)_{\infty}}\\\times
\frac{(a q/b_3b_{r+1},b_4,b_4/a,a q/b_3b_{r+2};q)_{\infty}}
{(b_3b_4/a,a^2q^2/b_3b_4b_{r+1}b_{r+2},q/b_{r+1},
a q/b_{r+1};q)_{\infty}}\\\times
\sum_{k=0} ^{\infty}\frac {(1-a b_3q^{1+2k}/b_4b_{r+1}b_{r+2})}
{(1-a b_3q/b_4b_{r+1}b_{r+2})}
\frac {(a b_3q/b_4b_{r+1}b_{r+2},a q/b_{r+1}b_{r+2},a q/b_4b_{r+1};q)_k}
{(q,b_3q/b_4,b_3q/b_{r+2};q)_k}\\\times
\frac {(a q/b_4b_{r+2},b_3,b_3/a;q)_k}
{(b_3q/b_{r+1},a q^2/b_4b_{r+1}b_{r+2},a^2 q^2/b_4b_{r+1}b_{r+2};q)_k}q^k\\
\times\Bigg[
\frac{(q,a,q/a,b_5,\dots,b_r,b_5/a,\dots,b_r/a;q)_{\infty}}
{(q^{1-k}/b_3,q/b_{r+3},\dots,q/b_{2r},
a q^{2+k}/b_4b_{r+1}b_{r+2};q)_{\infty}}\\\times
\frac{(b_3q^{1+k}/b_{r+3},\dots,b_3q^{1+k}/b_{2r},
a b_3q^{2+2k}/b_4b_{r+1}b_{r+2};q)_{\infty}}
{(a q^{1-k}/b_3,a q/b_{r+3},\dots,a q/b_{2r},
a^2 q^{2+k}/b_4b_{r+1}b_{r+2};q)_{\infty}}\\\times
\frac{(a q^{1-k}/b_3b_{r+3},\dots,
a q^{1-k}/b_3b_{2r},a^2 q^2/b_3b_4b_{r+1}b_{r+2};q)_{\infty}}
{(b_5q^{-k}/b_3,\dots,b_rq^{-k}/b_3,
b_3b_5q^k/a,\dots,b_3b_rq^k/a,
b_3^2q^{1+2k}/a;q)_{\infty}}\\\times
\sum_{n=0}^{\infty}\frac{(1-b_3^2q^{2k+2n}/a)}{(1-b_3^2q^{2k}/a)}
\frac{(b_3^2q^{2k}/a,b_3b_5q^k/a,\dots,b_3b_rq^k/a;q)_n}
{(q,b_3q^{1+k}/b_5,\dots,b_3q^{1+k}/b_r;q)_n}\\\times
\frac{(b_3b_{r+3}q^k/a,\dots,
b_3b_{2r}q^k/a,b_3b_4b_{r+1}b_{r+2}/a^2 q;q)_n}
{(b_3q^{1+k}/b_{r+3},\dots,
b_3q^{1+k}/b_{2r},a b_3q^{2+2k}/b_4b_{r+1}b_{r+2};q)_n}
\left(\frac{a^{r-1}q^{r-2}}{b_3\dots b_{2r}}\right)^n\\
+\Bigg(\frac{(q,a,q/a,b_3q^k,b_6,\dots,b_r,b_3q^k/a,
b_6/a,\dots,b_r/a;q)_{\infty}}
{(q/b_5,q/b_{r+3},\dots,q/b_{2r},
a q^{2+k}/b_4b_{r+1}b_{r+2};q)_{\infty}}\\\times
\frac{(b_5q/b_{r+3},\dots,b_5q/b_{2r},
a b_5q^{2+k}/b_4b_{r+1}b_{r+2};q)_{\infty}}
{(a q/b_5,a q/b_{r+3},\dots,a q/b_{2r},
a^2 q^{2+k}/b_4b_{r+1}b_{r+2};q)_{\infty}}\\\times
\frac{(a q/b_5b_{r+3},\dots,a q/b_5b_{2r},
a^2 q^{2+k}/b_4b_5b_{r+1}b_{r+2};q)_{\infty}}
{(b_3q^k/b_5,b_6/b_5,\dots,b_r/b_5,
b_3b_5q^k/a,b_5b_6/a,\dots,b_5b_r/a,
b_5^2q/a;q)_{\infty}}\\\times
\sum_{n=0}^{\infty}\frac{(1-b_5^2q^{2n}/a)}{(1-b_5^2/a)}
\frac{(b_5^2/a,b_3b_5q^k/a,b_5b_6/a,\dots,b_5b_r/a;q)_n}
{(q,b_5q^{1-k}/b_3,b_5q/b_6,\dots,b_5q/b_r;q)_n}\\\times
\frac{(b_5b_{r+3}/a,\dots,b_5b_{2r}/a,
b_4b_5b_{r+1}b_{r+2}q^{-1-k}/a^2;q)_n}
{(b_5q/b_{r+3},\dots,b_5q/b_{2r},a b_5q^{2+k}/b_4b_{r+1}b_{r+2};q)_n}
\left(\frac{a^{r-1}q^{r-2}}{b_3\dots b_{2r}}\right)^n\\
+\operatorname{idem}(b_5;b_6,\dots,b_r)\Bigg)\Bigg]
+\operatorname{idem}(b_3;b_4),
\end{multline*}
which can be simplified to
\begin{multline}\label{2r2rgl3}
\Bigg[\frac{(b_3q/b_{r+2},b_3q/b_{r+1},a q/b_3b_{r+1},b_4,b_4/a,
a q/b_3b_{r+2},q,a,q/a,b_5,\dots,b_r;q)_{\infty}}
{(b_4/b_3,q/b_{r+2},
a q/b_{r+2},b_3b_4/a,q/b_{r+1},a q/b_{r+1},
q/b_3,q/b_{r+3},\dots,q/b_{2r};q)_{\infty}}\\\times
\frac{(b_5/a,\dots,b_r/a,b_3q/b_{r+3},\dots,b_3q/b_{2r},
a q/b_3b_{r+3},\dots,a q/b_3b_{2r};q)_{\infty}}
{(a q/b_3,a q/b_{r+3},\dots,a q/b_{2r},
b_5/b_3,\dots,b_r/b_3,b_3b_5/a,\dots,b_3b_r/a,
b_3^2q/a;q)_{\infty}}\\\times
\sum_{k=0} ^{\infty}\frac {(1-a b_3q^{1+2k}/b_4b_{r+1}b_{r+2})}
{(1-a b_3q/b_4b_{r+1}b_{r+2})}
\frac {(a b_3q/b_4b_{r+1}b_{r+2},a q/b_{r+1}b_{r+2},a q/b_4b_{r+1};q)_k}
{(q,b_3q/b_4,b_3q/b_{r+2};q)_k}\\\times
\frac {(a q/b_4b_{r+2},b_3b_5/a,\dots,b_3b_r/a,
b_3b_{r+3}/a,\dots,b_3b_{2r}/a;q)_k}
{(b_3q/b_{r+1},b_3q/b_5,\dots,b_3q/b_r,b_3q/b_{r+3},\dots,
b_3q/b_{2r};q)_k}\\\times
\frac {(b_3^2q/a;q)_{2k}}
{(a b_3q^2/b_4b_{r+1}b_{r+2};q)_{2k}}
\left(\frac{a^{r-3}q^{r-3}}{b_5\dots b_rb_{r+3}\dots b_{2r}}\right)^k\\\times
\sum_{n=0}^{\infty}\frac{(1-b_3^2q^{2k+2n}/a)}{(1-b_3^2q^{2k}/a)}
\frac{(b_3^2q^{2k}/a,b_3b_5q^k/a,\dots,b_3b_rq^k/a;q)_n}
{(q,b_3q^{1+k}/b_5,\dots,b_3q^{1+k}/b_r;q)_n}\\\times
\frac{(b_3b_{r+3}q^k/a,\dots,
b_3b_{2r}q^k/a,b_3b_4b_{r+1}b_{r+2}/a^2 q;q)_n}
{(b_3q^{1+k}/b_{r+3},\dots,
b_3q^{1+k}/b_{2r},a b_3q^{2+2k}/b_4b_{r+1}b_{r+2};q)_n}
\left(\frac{a^{r-1}q^{r-2}}{b_3\dots b_{2r}}\right)^n\\
+\Bigg(
\frac{(b_3q/b_{r+2},b_3q/b_{r+1},a q/b_3b_{r+1},b_4,b_4/a,
a q/b_3b_{r+2},q,a,q/a;q)_{\infty}}
{(a b_3q^2/b_4b_{r+1}b_{r+2},b_4/b_3,q/b_{r+2},
a q/b_{r+2},b_3b_4/a,a^2q^2/b_3b_4b_{r+1}b_{r+2};q)_{\infty}}\\\times
\frac{(b_3,b_6,\dots,b_r,b_3/a,
b_6/a,\dots,b_r/a,b_5q/b_{r+3},\dots,b_5q/b_{2r};q)_{\infty}}
{(q/b_{r+1},a q/b_{r+1},q/b_5,q/b_{r+3},\dots,q/b_{2r},
a q/b_5,a q/b_{r+3},\dots,a q/b_{2r};q)_{\infty}}\\\times
\frac{(a b_5q^2/b_4b_{r+1}b_{r+2},a q/b_5b_{r+3},\dots,a q/b_5b_{2r},
a^2 q^2/b_4b_5b_{r+1}b_{r+2};q)_{\infty}}
{(b_3/b_5,b_6/b_5,\dots,b_r/b_5,b_3b_5/a,b_5b_6/a,\dots,b_5b_r/a,
b_5^2q/a;q)_{\infty}}\\\times
\sum_{k=0} ^{\infty}\frac {(1-a b_3q^{1+2k}/b_4b_{r+1}b_{r+2})}
{(1-a b_3q/b_4b_{r+1}b_{r+2})}
\frac {(a b_3q/b_4b_{r+1}b_{r+2},a q/b_{r+1}b_{r+2},a q/b_4b_{r+1};q)_k}
{(q,b_3q/b_4,b_3q/b_{r+2};q)_k}\\\times
\frac {(a q/b_4b_{r+2},b_3/b_5,b_3b_5/a;q)_k}
{(b_3q/b_{r+1},a b_5q^2/b_4b_{r+1}b_{r+2},
a^2 q^2/b_4b_5b_{r+1}b_{r+2};q)_k}q^k\\\times
\sum_{n=0}^{\infty}\frac{(1-b_5^2q^{2n}/a)}{(1-b_5^2/a)}
\frac{(b_5^2/a,b_3b_5q^k/a,b_5b_6/a,\dots,b_5b_r/a;q)_n}
{(q,b_5q^{1-k}/b_3,b_5q/b_6,\dots,b_5q/b_r;q)_n}\\\times
\frac{(b_5b_{r+3}/a,\dots,b_5b_{2r}/a,
b_4b_5b_{r+1}b_{r+2}q^{-1-k}/a^2;q)_n}
{(b_5q/b_{r+3},\dots,b_5q/b_{2r},a b_5q^{2+k}/b_4b_{r+1}b_{r+2};q)_n}
\left(\frac{a^{r-1}q^{r-2}}{b_3\dots b_{2r}}\right)^n\\
+\operatorname{idem}(b_5;b_6,\dots,b_r)\Bigg)\Bigg]
+\operatorname{idem}(b_3;b_4).
\end{multline}

We have in \eqref{2r2rgl3} a sum of $2(r-3)$ double sums.
Accordingly, for more clarity, let us write the whole expression
in \eqref{2r2rgl3} as
\begin{align}\label{tu}\notag
T_1(b_3,b_4)+{} &T_2(b_3,b_4)+\dots +T_{r-3}(b_3,b_4)\\
{}+ { } U_1(b_3,b_4)+{} &U_2(b_3,b_4)+\dots +U_{r-3}(b_3,b_4),
\end{align}
where, by definition of ``idem", $U_i(b_3,b_4)=T_i(b_4,b_3)$
for $i=1,\dots,r-3$. Further, $\sum_{i=2}^{r-3}T_i(b_3,b_4)
=T_2(b_3,b_4)+\operatorname{idem}(b_5;b_6,\dots,b_r)$
(and $\sum_{i=2}^{r-3}U_i(b_3,b_4)
=U_2(b_3,b_4)+\operatorname{idem}(b_5;b_6,\dots,b_r)$).
Below, we will selectively perform manipulations with the respective terms
$T_i$ and $U_i$ ($i=1,\dots,r-3$). 

To evaluate $T_1(b_3,b_4)$ (and hence also $U_1(b_3,b_4)$),
we first shift the index $n$ of the inner sum in $T_1(b_3,b_4)$ (or,
equivalently, in the first term in \eqref{2r2rgl3}) by $-k$
and then interchange the double sum.
Symbolically, we apply the interchange of summations as in \eqref{intchs}.
Thus, we obtain (using some elementary identities for $q$-shifted factorials)
\begin{multline*}
T_1(b_3,b_4)
=\frac{(q,a,q/a,b_4,\dots,b_r,b_4/a,\dots,b_r/a;q)_{\infty}}
{(q/b_3,q/b_{r+1},\dots,q/b_{2r},a q/b_3,a q/b_{r+1},\dots,
a q/b_{2r};q)_{\infty}}\\\times
\frac{(b_3q/b_{r+1},\dots,b_3q/b_{2r},
a q/b_3b_{r+1},\dots,a q/b_3b_{2r};q)_{\infty}}
{(b_4/b_3,\dots,b_r/b_3,b_3b_4/a,\dots,b_3b_r/a,
b_3^2q/a;q)_{\infty}}\\\times
\sum_{n=0}^{\infty}\frac{(1-b_3^2q^{2n}/a)}{(1-b_3^2/a)}
\frac{(b_3^2/a,b_3b_5/a,\dots,b_3b_r/a;q)_n}
{(q,b_3q/b_5,\dots,b_3q/b_r;q)_n}\\\times
\frac{(b_3b_{r+3}/a,\dots,
b_3b_{2r}/a,b_3b_4b_{r+1}b_{r+2}/a^2 q;q)_n}
{(b_3q/b_{r+3},\dots,
b_3q/b_{2r},a b_3q^2/b_4b_{r+1}b_{r+2};q)_n}
\left(\frac{a^{r-1}q^{r-2}}{b_3\dots b_{2r}}\right)^n\\\times
\sum_{k=0} ^n\frac {(1-a b_3q^{1+2k}/b_4b_{r+1}b_{r+2})}
{(1-a b_3q/b_4b_{r+1}b_{r+2})}
\frac {(a b_3q/b_4b_{r+1}b_{r+2},a q/b_{r+1}b_{r+2},a q/b_4b_{r+1};q)_k}
{(q,b_3q/b_4,b_3q/b_{r+2};q)_k}\\\times
\frac {(a q/b_4b_{r+2},b_3^2q^n/a,q^{-n};q)_k}
{(b_3q/b_{r+1},a^2q^{2-n}/b_3b_4b_{r+1}b_{r+2},
a b_3q^{2+n}/b_4b_{r+1}b_{r+2};q)_k}q^k.
\end{multline*}
Now the inner sum can be evaluated by Jackson's terminating
$_8\phi_7$ summation in \eqref{87gl}, which simplifies the last expression,
$T_1(b_3,b_4)$, to
\begin{multline}\label{2r2rgl4}
\frac{(q,a,q/a,b_4,\dots,b_r,b_4/a,\dots,b_r/a;q)_{\infty}}
{(q/b_3,q/b_{r+1},\dots,q/b_{2r},a q/b_3,a q/b_{r+1},\dots,
a q/b_{2r};q)_{\infty}}\\\times
\frac{(b_3q/b_{r+1},\dots,b_3q/b_{2r},
a q/b_3b_{r+1},\dots,a q/b_3b_{2r};q)_{\infty}}
{(b_4/b_3,\dots,b_r/b_3,b_3b_4/a,\dots,b_3b_r/a,
b_3^2q/a;q)_{\infty}}\\\times
\sum_{n=0}^{\infty}\frac{(1-b_3^2q^{2n}/a)}{(1-b_3^2/a)}
\frac{(b_3^2/a,b_3b_4/a,\dots,b_3b_{2r}/a;q)_n}
{(q,b_3q/b_4,\dots,b_3q/b_{2r};q)_n}
\left(\frac{a^{r-1}q^{r-2}}{b_3\dots b_{2r}}\right)^n.
\end{multline}

Next, we consider $T_2(b_3,b_4)$ and $U_2(b_3,b_4)$.
By interchanging the double sums in $T_2(b_3,b_4)$ and in $U_2(b_3,b_4)$
we obtain
\begin{multline}\label{2r2rgl5}
T_2(b_3,b_4) + U_2(b_3,b_4)\\
=\frac{(b_3q/b_{r+2},b_3q/b_{r+1},a q/b_3b_{r+1},b_4,b_4/a,
a q/b_3b_{r+2},q,a,q/a;q)_{\infty}}
{(a b_3q^2/b_4b_{r+1}b_{r+2},b_4/b_3,q/b_{r+2},
a q/b_{r+2},b_3b_4/a,a^2q^2/b_3b_4b_{r+1}b_{r+2};q)_{\infty}}\\\times
\frac{(b_3,b_6,\dots,b_r,b_3/a,
b_6/a,\dots,b_r/a,b_5q/b_{r+3},\dots,b_5q/b_{2r};q)_{\infty}}
{(q/b_{r+1},a q/b_{r+1},q/b_5,q/b_{r+3},\dots,q/b_{2r},
a q/b_5,a q/b_{r+3},\dots,a q/b_{2r};q)_{\infty}}\\\times
\frac{(a b_5q^2/b_4b_{r+1}b_{r+2},a q/b_5b_{r+3},\dots,a q/b_5b_{2r},
a^2 q^2/b_4b_5b_{r+1}b_{r+2};q)_{\infty}}
{(b_3/b_5,b_6/b_5,\dots,b_r/b_5,b_3b_5/a,b_5b_6/a,\dots,b_5b_r/a,
b_5^2q/a;q)_{\infty}}\\\times
\sum_{n=0}^{\infty}\frac{(1-b_5^2q^{2n}/a)}{(1-b_5^2/a)}
\frac{(b_5^2/a,b_3b_5/a,b_5b_6/a,\dots,b_5b_r/a;q)_n}
{(q,b_5q/b_3,b_5q/b_6,\dots,b_5q/b_r;q)_n}\\\times
\frac{(b_5b_{r+3}/a,\dots,b_5b_{2r}/a,
b_4b_5b_{r+1}b_{r+2}/a^2q;q)_n}
{(b_5q/b_{r+3},\dots,b_5q/b_{2r},a b_5q^2/b_4b_{r+1}b_{r+2};q)_n}
\left(\frac{a^{r-1}q^{r-2}}{b_3\dots b_{2r}}\right)^n\\\times
\sum_{k=0} ^{\infty}\frac {(1-a b_3q^{1+2k}/b_4b_{r+1}b_{r+2})}
{(1-a b_3q/b_4b_{r+1}b_{r+2})}
\frac {(a b_3q/b_4b_{r+1}b_{r+2},a q/b_{r+1}b_{r+2},a q/b_4b_{r+1};q)_k}
{(q,b_3q/b_4,b_3q/b_{r+2};q)_k}\\\times
\frac {(a q/b_4b_{r+2},b_3q^{-n}/b_5,b_3b_5q^n/a;q)_k}
{(b_3q/b_{r+1},a b_5q^{2+n}/b_4b_{r+1}b_{r+2},
a^2 q^{2-n}/b_4b_5b_{r+1}b_{r+2};q)_k}q^k\\
+\operatorname{idem}(b_3;b_4).
\end{multline}
Now, to the inner sum of the {\em first} double sum (but not of the second!)
in \eqref{2r2rgl5} we apply
Bailey's nonterminating $_8\phi_7$ summation \eqref{87ntgl}, i.e.,
specifically we apply
\begin{multline}\label{87ntgl2r}
\sum_{k=0}^{\infty}\frac{(1-a b_3q^{1+2k}/b_4b_{r+1}b_{r+2})}
{((1-a b_3q/b_4b_{r+1}b_{r+2}))}
\frac{(a b_3q/b_4b_{r+1}b_{r+2},a q/b_{r+1}b_{r+2};q)_k}
{(q,b_3q/b_4,b_3q/b_{r+2};q)_k}\\\times
\frac{(a q/b_4b_{r+1},
a q/b_4b_{r+2},b_3q^{-n}/b_5,b_3b_5q^n/a;q)_k}
{(b_3q/b_{r+1},a b_5q^{2+n}/b_4b_{r+1}b_{r+2},
a^2q^{2-n}/b_4b_5b_{r+1}b_{r+2};q)_k}q^k\\
=\frac{(a b_3q^2/b_4b_{r+1}b_{r+2},b_4/b_3,b_3b_4/a,
b_5q^{1+n}/b_{r+2};q)_{\infty}}
{(b_3q/b_{r+2},b_3q/b_{r+1},a b_5q^{2+n}/b_4b_{r+1}b_{r+2},
a^2q^{2-n}/b_4b_5b_{r+1}b_{r+2};q)_{\infty}}\\\times
\frac{(a q^{1-n}/b_5b_{r+2},b_5q^{1+n}/b_{r+1},a q^{1-n}/b_5b_{r+1},
a^2q^2/b_3b_4b_{r+1}b_{r+2};q)_{\infty}}
{(a q/b_3b_{r+1},a q/b_3b_{r+2},b_4q^{-n}/b_5,
b_4b_5q^n/a;q)_{\infty}}\\
-\frac{(a b_3q^2/b_4b_{r+1}b_{r+2},a q/b_4b_{r+1},
a q/b_4b_{r+2},b_3q^{-n}/b_5,b_3b_5q^n/a;q)_{\infty}}
{(b_3/b_4,b_3q/b_{r+2},b_3q/b_{r+1},a b_5q^{2+n}/b_4b_{r+1}b_{r+2},
a ^2q^{2-n}/b_4b_5b_{r+1}b_{r+2};q)_{\infty}}\\\times
\frac{(b_4/b_3,b_4q/b_{r+2},b_4q/b_{r+1},
a b_5q^{2+n}/b_3b_{r+1}b_{r+2},a ^2q^{2-n}/b_3b_5b_{r+1}b_{r+2};q)_{\infty}}
{(a q/b_3b_{r+1},a q/b_3b_{r+2},b_4q^{-n}/b_5,
b_4b_5q^n/a,a b_4q^2/b_3b_{r+1}b_{r+2};q)_{\infty}}\\\times
\sum_{k=0}^{\infty}\frac{(1-a b_4q^{1+2k}/b_3b_{r+1}b_{r+2})}
{(1-a b_4q/b_3b_{r+1}b_{r+2})}
\frac{(a b_4q/b_3b_{r+1}b_{r+2},a q/b_{r+1}b_{r+2};q)_k}
{(q,b_4q/b_3;q)_k}\\\times
\frac{(a q/b_3b_{r+1},
a q/b_3b_{r+2},b_4q^{-n}/b_5,b_4b_5q^n/a ;q)_k}
{(b_4q/b_{r+2},b_4q/b_{r+1},a b_5q^{2+n}/b_3b_{r+1}b_{r+2},
a ^2q^{2-n}/b_3b_5b_{r+1}b_{r+2};q)_k}q^k.
\end{multline}
The result of the application of this (two term) summation is that
the first term in \eqref{2r2rgl5}, $T_2(b_3,b_4)$, is split into
two parts, one single sum $T_2^{\prime}(b_3,b_4)$ and one double
sum $T_2^{\prime\prime}(b_3,b_4)$.
Formally, we have
$$
T_2(b_3,b_4)+U_2(b_3,b_4)=[T_2^{\prime}(b_3,b_4)+
T_2^{\prime\prime}(b_3,b_4)]+U_2(b_3,b_4).
$$
But $T_2^{\prime\prime}(b_3,b_4)$ is precisely $-U_2(b_3,b_4)$
(as can be readily checked), so two terms cancel, thus
$$
T_2(b_3,b_4)+U_2(b_3,b_4)=T_2^{\prime}(b_3,b_4).
$$
Now, the last expression, $T_2^{\prime}(b_3,b_4)$, can be simplified to
\begin{multline}\label{2r2rgl6}
\frac{(q,a,q/a,b_3,b_4,b_6,\dots,b_r,b_3/a,b_4/a,
b_6/a,\dots,b_r/a;q)_{\infty}}
{(q/b_5,q/b_{r+1},\dots,q/b_{2r},
a q/b_5,a q/b_{r+1},\dots,a q/b_{2r},b_3/b_5,b_4/b_5;q)_{\infty}}\\\times
\frac{(b_5q/b_{r+1},\dots,b_5q/b_{2r},
a q/b_5b_{r+1},\dots,a q/b_5b_{2r};q)_{\infty}}
{(b_6/b_5,\dots,b_r/b_5,b_3b_5/a,b_4b_5/a,
b_5b_6/a,\dots,b_5b_r/a,b_5^2q/a;q)_{\infty}}\\\times
\sum_{n=0}^{\infty}\frac{(1-b_5^2q^{2n}/a)}{(1-b_5^2/a)}
\frac{(b_5^2/a,b_3b_5/a,b_4b_5/a;q)_n}
{(q,b_5q/b_3,b_5q/b_4;q)_n}\\\times
\frac{(b_5b_6/a,\dots,b_5b_{2r}/a;q)_n}
{(b_5q/b_6,\dots,b_5q/b_{2r};q)_n}
\left(\frac{a^{r-1}q^{r-2}}{b_3\dots b_{2r}}\right)^n.
\end{multline}
It is easy to see that \eqref{2r2rgl6} equals
\eqref{2r2rgl4} where $b_3$ and $b_5$ are interchanged.
Collecting all terms, according to \eqref{tu}, completes our derivation
of the very-well-poised $_{2r}\psi_{2r}$ transformation formula
in \eqref{2r2rgl}.

\medskip
We conclude this section considering Slater's well-poised $_{2r}\psi_{2r}$
transformation in its general form.

If, in \eqref{2r2rgl}, we set $b_{2r-1}=-b_{2r}=\sqrt{a}$, and
then let $r\mapsto r+2$, and relabel $b_i\mapsto b_{i-2}$, we
obtain the following transformation formula for
a well-poised $_{2r}\psi_{2r}$ series:
\begin{multline}\label{slw2r2rgl}
{}_{2r}\psi_{2r}\!\left[\begin{matrix}b_1,b_2,
\dots,b_{2r}\\
\frac{a q}{b_1},\frac{a q}{b_2},
\dots,\frac{a q}{b_{2r}}\end{matrix}\,;q,
-\frac{a^rq^r}{b_1\dots b_{2r}}\right]\\
=\frac{(q,a,\frac q{a},b_2,\dots,b_r,\frac{b_2}{a},\dots,\frac{b_r}{a},
\frac{b_1q}{b_{r+1}},\dots,\frac{b_1q}{b_{2r}},
\frac{a q}{b_1b_{r+1}},\dots,\frac{a q}{b_1b_{2r}},;q)_{\infty}}
{(\frac q{b_1},\frac q{b_{r+1}},\dots,\frac q{b_{2r}},
\frac{a q}{b_1},\frac{a q}{b_{r+1}},\dots,\frac{a q}{b_{2r}},
\frac{b_2}{b_1},\dots,\frac{b_r}{b_1},\frac{b_1b_2}{a},\dots,
\frac{b_1b_r}{a},\frac{b_1^2q}{a};q)_{\infty}}\\\times
\frac{(\frac{b_1q}{\sqrt{a}},-\frac{b_1q}{\sqrt{a}},\frac{\sqrt{a}q}{b_1},
-\frac{\sqrt{a}q}{b_1};q)_{\infty}}
{(\frac q{\sqrt{a}},-\frac q{\sqrt{a}},\sqrt{a},
-\sqrt{a};q)_{\infty}}\;
{}_{2r}\phi_{2r-1}\!\left[\begin{matrix}\frac{b_1^2}{a},
\frac{b_1b_2}{a},\frac{b_1b_3}{a},\dots,\frac{b_1b_{2r}}{a}\\
\frac{b_1q}{b_2},\frac{b_1q}{b_3},\dots,\frac{b_1q}{b_{2r}}\end{matrix}\,;q,
-\frac{a^rq^r}{b_1\dots b_{2r}}\right]\\
+\operatorname{idem}(b_1;b_2,\dots,b_r),
\end{multline}
where the series either terminate, or $|a^rq^r/b_1\dots b_{2r}|<1$,
for convergence. On the other hand, we can derive \eqref{2r2rgl}
from \eqref{slw2r2rgl} by choosing $b_{2r-1}=-b_{2r}=q\sqrt{a}$
and relabelling of the parameters $b_i\mapsto b_{i+2}$.
Thus, the transformations \eqref{2r2rgl} and \eqref{slw2r2rgl} are equivalent.

Slater~\cite{slatertf}, in fact, derived the more general transformation
\begin{multline}\label{sl2r2rgl}
{}_{2r}\psi_{2r}\!\left[\begin{matrix}b_1,b_2,
\dots,b_{2r}\\
\frac{a q}{b_1},\frac{a q}{b_2},
\dots,\frac{a q}{b_{2r}}\end{matrix}\,;q,
-\frac{a^rq^r}{b_1\dots b_{2r}}\right]\\
=\frac{(a,\frac q{a},a_2,\dots,a_r,\frac q{a_2},\dots,\frac q{a_r},
\frac{a_2}{a},\dots,\frac{a_r}{a},
\frac{aq}{a_2},\dots,\frac{aq}{a_r};q)_{\infty}}
{(\frac q{b_1},\dots,\frac q{b_{2r}},
\frac{a q}{b_1},\dots,\frac{a q}{b_{2r}},
\frac{a_2}{a_1},\dots,\frac{a_r}{a_1},\frac{a_1q}{a_2},\dots,
\frac{a_1q}{a_r};q)_{\infty}}\\\times
\frac{(\frac{a_1q}{b_1},\dots,\frac{a_1q}{b_{2r}},
\frac{a q}{a_1b_1},\dots,\frac{a q}{a_1b_{2r}},
\frac{a_1}{\sqrt{a}},-\frac{a_1}{\sqrt{a}},\frac{\sqrt{a}q}{a_1},
-\frac{\sqrt{a}q}{a_1};q)_{\infty}}
{(\frac{a_1a_2}{a},\dots,
\frac{a_1a_r}{a},\frac{aq}{a_1a_2},\dots,\frac{aq}{a_1a_r},
\frac{a_1^2}{a},\frac{aq}{a_1^2},
\frac q{\sqrt{a}},-\frac q{\sqrt{a}},\sqrt{a},-\sqrt{a};q)_{\infty}}\\\times
{}_{2r}\psi_{2r}\!\left[\begin{matrix}\frac{a_1b_1}{a},
\frac{a_1b_2}{a},\dots,\frac{a_1b_{2r}}{a}\\
\frac{a_1q}{b_1},\frac{a_1q}{b_2},\dots,\frac{a_1q}{b_{2r}}\end{matrix}\,;q,
-\frac{a^rq^r}{b_1\dots b_{2r}}\right]
+\operatorname{idem}(a_1;a_2,\dots,a_r),
\end{multline}
where the series either terminate, or $|a^rq^r/b_1\dots b_{2r}|<1$,
for convergence.

It is easy to see, that the transformation in \eqref{sl2r2rgl}
which involves only bilateral series is much more general than
the transformation in \eqref{slw2r2rgl} since in \eqref{sl2r2rgl}
we have $r$ additional parameters $a_1,\dots,a_r$. The special
case $a_i = b_i$, $i=1,\dots,r$, of \eqref{sl2r2rgl} is
exactly \eqref{slw2r2rgl}.

We were able to prove the transformation in \eqref{slw2r2rgl}
(which is equivalent to \eqref{2r2rgl}) by induction. We started from
less complicated identities and ultimatively deduced more complex ones
involving more parameters.
The natural question arises whether we can also prove the more general
transformation \eqref{sl2r2rgl} by elementary means.
Unfortunately, we were not able to derive \eqref{sl2r2rgl} directly by
the method of this article. Instead, we want to point out that
the transformation we derived, \eqref{slw2r2rgl}, can be extended
to \eqref{sl2r2rgl} by an $r$-fold application of ``Ismail's~\cite{ismail}
argument" (which is actually a classical analytic continuation argument).
We provide a sketch of how this works.

It is easy to see that both sides of \eqref{sl2r2rgl} are analytic in
each $b_1^{-1},b_2^{-1},\dots,b_r^{-1}$ in a domain around the origin.
We know that the identity is true when $b_i=a_i$, for $i=1,\dots,r$.
What needs to be done is to extend \eqref{slw2r2rgl} first by an additional
parameter $b_1$, then by $b_2$, etc. This means that if we have
already extended \eqref{slw2r2rgl} by $b_1,\dots,b_j$, what we
should have derived is the identity in \eqref{sl2r2rgl} where $b_i=a_i$,
for $i=j+1,\dots,r$. So, we proceed by induction starting with $j=0$
(identity \eqref{slw2r2rgl}) and ending with $j=r$ (identity
\eqref{sl2r2rgl}). In the inductive step, we consider the
transformation \eqref{sl2r2rgl} where $b_i=a_i$, for $i=j+1,\dots,r$.
We call that identity $E_{j+1}$. We need to show that $E_{j+1}$ is true,
provided $E_j$ is true, which is \eqref{sl2r2rgl} where $b_i=a_i$,
for $i=j,\dots,r$.
We observe that both sides of $E_{j+1}$ are
analytic in $1/b_j$ around the origin.
The identity is true for $b_j=a_jq^{-n_j}$, $n_j=0,1,2,\dots$.
This follows by the $a_j\mapsto a_jq^{-n_j}$ case
of the inductive hypothesis, $E_j$. (This can be verified by looking at all
the terms involving $j$. Further note that the index of the
$j$-th sum is shifted by $-n_j$, since the $j$-th bilateral series
becomes a unilateral series.) Since $E_{j+1}$ is true for
an infinite sequence of $1/b_j$ which has an accumulation point,
namely 0, in the interior of the domain $\mathcal D$ of analyticity of
$1/b_j$, we can apply the identity theorem to deduce that $E_{j+1}$
is true for $1/b_j$ throughout $\mathcal D$.
Now, by induction, \eqref{sl2r2rgl} follows, with the general additional
parameters $b_1,\dots,b_r$.

\section{Slater's general bilateral transformations}\label{secslgen}

Here we consider series which are not necessarily well-poised.
An important transformation for general $_r\psi_r$ series was given by
Slater~\cite{slatertf}. Her formula in its general form,
see \eqref{slgentf}, connects $r+1$ bilateral $_r\psi_r$ series.
However, the special case of the transformation
where only one of the series is bilateral
and the $r$ other series are all unilateral ($_r\phi_{r-1}$ series)
itself is interesting, as it includes many important identities
for (bilateral) basic hypergeometric series as special cases. 
This formula reads as follows:
\begin{multline}\label{slgentfu}
{}_r\psi_r\!\left[\begin{matrix}a_1,a_2,\dots,a_r\\b_1,b_2,
\dots,b_r\end{matrix}\,;q,z\right]\\
=\frac{(q,a_2,\dots,a_r,\frac{b_1}{a_1},\dots,\frac{b_r}{a_1},
a_1z,\frac q{a_1z};q)_{\infty}}
{(\frac q{a_1},\frac{a_2}{a_1},\dots,\frac{a_r}{a_1},
b_1,\dots,b_r,z,\frac qz;q)_{\infty}}\;
{}_r\phi_{r-1}\!\left[\begin{matrix}\frac {a_1q}{b_1},\frac{a_1q}{b_2},
\dots,\frac {a_1q}{b_r}\\
\frac {a_1q}{a_2},\dots,\frac{a_1q}{a_r}\end{matrix}\,;q,\frac{b_1\dots b_r}
{a_1\dots a_rz}\right]\\
+\operatorname{idem}(a_1;a_2,\dots,a_r),
\end{multline}
where the series either terminate, or $|b_1\dots b_r/a_1\dots a_r|<|z|<1$,
for convergence.
(The symbol ``$\operatorname{idem}(a_1;a_2,\dots,a_r)$" is explained in
the introduction.)

We are able to give an elementary inductive proof of this general
transformation formula. We proceed by similar means as in the previous
sections where Slater's general transformation for
well-poised bilateral $_{2r}\psi_{2r}$ series was derived.
For deriving the transformation in \eqref{slgentfu} we make use of other
identities. In particular, where we were before using identities for
very-well-poised $_8\phi_7$ series,
here we instead utilize identities for balanced $_3\phi_2$ series.

The $r=1$ case of \eqref{slgentfu} is readily verifed
using the $_1\psi_1$ and $_1\phi_0$ summations in
\eqref{11gl} and in \eqref{10gl}, respectively.
Ramanujan's~\cite{hardy} $_1\psi_1$ summation
(cf.~\cite[Eq.~(5.2.1)]{grhyp}) reads
\begin{equation}\label{11gl}
{}_1\psi_1\!\left[\begin{matrix}a\\
b\end{matrix};q,z\right]=
\frac{(q,b/a,az,q/az;q)_{\infty}}{(b,q/a,z,b/az;q)_{\infty}},
\end{equation}
where the series either terminates, or $|b/a|<|z|<1$, for convergence.
The summation in \eqref{11gl} is a bilateral extension of the $q$-binomial
theorem (cf.~\cite[Sec.~1.3]{grhyp}),
\begin{equation}\label{10gl}
{}_1\phi_0\!\left[\begin{matrix}a\\
-\end{matrix};q,z\right]=
\frac{(az;q)_{\infty}}{(z;q)_{\infty}},
\end{equation}
where the series either terminates, or $|z|<1$, for convergence.
The summation in \eqref{10gl} was first discovered by Cauchy~\cite{cauchy}.
Clearly, \eqref{11gl} reduces to \eqref{10gl} when $b=q$.
The first elementary proof of the $_1\psi_1$ summation formula~\eqref{11gl}
was given by M.~Jackson (as pointed out to us by George
Andrews~\cite{andpriv}). Jackson's proof essentially derives the
$_1\psi_1$ summation from the $q$-Gau{\ss} summation, by manipulation
of series. In \cite{schlelsum} we reviewed Jackson's proof and extended it
to a method to provide new elementary proofs of Dougall's~\cite{dougall}
$_2H_2$ summation and Bailey's~\cite{bail66} very-well-poised
$_6\psi_6$ summation, respectively. 

Now let us establish the $r=2$ case of \eqref{slgentfu}.

Here, and later we make use of a special case of the identity in
\eqref{87ntgl2}. Namely set $b\mapsto a_2$, $c\mapsto aq/b_2$,
$d\mapsto a_1$, and $e\mapsto aq/b_1$ in
\eqref{87ntgl2}, and then let $a\to\infty$. We obtain
\begin{multline}\label{key32gl}
\frac{(b_1b_2/a_2,b_2/a_1,a_2,b_1/a_1;q)_{\infty}}
{(a_2/a_1,b_1,b_1b_2/a_1a_2,b_2;q)_{\infty}}
\sum_{k=0}^{\infty}\frac{(b_1/a_2,b_2/a_2;q)_k(a_1;q)_{n+k}}
{(q,a_1q/a_2;q)_k(b_1b_2/a_2;q)_{n+k}}q^k\\
+\operatorname{idem}(a_1;a_2)=
\frac{(a_1,a_2;q)_n}{(b_1,b_2;q)_n}.
\end{multline}
Alternatively, we could also have invoked the non-terminating
$q$-Pfaff--Saalsch\"utz summation theorem (cf.~\cite[Eq.~(II.24)]{grhyp}),
\begin{multline}\label{32ntgl}
{}_3\phi_2\!\left[\begin{matrix}
a,b,c\\e,f\end{matrix};\,q,q\right]+
\frac{(q/e,a,b,c,fq/e;q)_{\infty}}{(e/q,aq/e,bq/e,cq/e,f;q)_{\infty}}\,
{}_3\phi_2\!\left[\begin{matrix}
aq/e,bq/e,cq/e\\q^2/e,fq/e\end{matrix};\,q,q\right]\\
=\frac{(q/e,f/a,f/b,f/c;q)_{\infty}}{(aq/e,bq/e,cq/e,f;q)_{\infty}},
\end{multline}
where $ef=abcq$. If we simultaneously replace $a$, $b$, $c$,
and $e$ in \eqref{32ntgl} by $b_1/a_2$, $b_2/a_2$, $a_1q^n$, and
$a_1q/a_2$, respectively,
we obtain after some manipulations \eqref{key32gl}.

For the $r=2$ case of \eqref{slgentfu}, we multiply both sides of
\eqref{key32gl} by $z^n$ and sum over all integers $n$.
On the right side we obtain
\begin{equation*}
{}_{2}\psi_{2}\!\left[\begin{matrix}a_1,a_2\\b_1,b_2\end{matrix}\,;q,
z\right].
\end{equation*}
On the left side we obtain
\begin{multline*}
\frac{(b_1b_2/a_2,b_2/a_1,a_2,b_1/a_1;q)_{\infty}}
{(a_2/a_1,b_1,b_1b_2/a_1a_2,b_2;q)_{\infty}}
\sum_{n=-\infty}^{\infty}z^n
\sum_{k=0}^{\infty}\frac{(b_1/a_2,b_2/a_2;q)_k(a_1;q)_{n+k}}
{(q,a_1q/a_2;q)_k(b_1b_2/a_2;q)_{n+k}}q^k\\
+\operatorname{idem}(a_1;a_2)=
\frac{(b_1b_2/a_2,b_2/a_1,a_2,b_1/a_1;q)_{\infty}}
{(a_2/a_1,b_1,b_1b_2/a_1a_2,b_2;q)_{\infty}}
\sum_{k=0}^{\infty}\frac{(b_1/a_2,b_2/a_2,a_1;q)_k}
{(q,a_1q/a_2,b_1b_2/a_2;q)_k}q^k\\\times
\sum_{n=-\infty}^{\infty}\frac{(a_1q^k;q)_n}
{(b_1b_2q^k/a_2;q)_n}z^n
+\operatorname{idem}(a_1;a_2).
\end{multline*}
Now we can evaluate the the inner sums by the $_1\psi_1$ summation~\eqref{11gl}.
This yields
\begin{multline*}
\frac{(b_1b_2/a_2,b_2/a_1,a_2,b_1/a_1;q)_{\infty}}
{(a_2/a_1,b_1,b_1b_2/a_1a_2,b_2;q)_{\infty}}
\sum_{k=0}^{\infty}\frac{(b_1/a_2,b_2/a_2,a_1;q)_k}
{(q,a_1q/a_2,b_1b_2/a_2;q)_k}q^k\\\times
\frac{(q,b_1b_2/a_1a_2,a_1zq^k,q^{1-k}/a_1z;q)_{\infty}}
{(b_1b_2q^k/a_2,q^{1-k}/a_1,z,b_1b_2/a_1a_2z;q)_{\infty}}
+\operatorname{idem}(a_1;a_2),
\end{multline*}
which can be simplified to
\begin{multline}\label{2psi2gl}
\frac{(b_2/a_1,a_2,b_1/a_1,q,a_1z,q/a_1z;q)_{\infty}}
{(a_2/a_1,b_1,b_2,q/a_1,z,b_1b_2/a_1a_2z;q)_{\infty}}
\sum_{k=0}^{\infty}\frac{(b_1/a_2,b_2/a_2;q)_k}
{(q,a_1q/a_2;q)_k}\left(\frac qz\right)^k\\
+\operatorname{idem}(a_1;a_2).
\end{multline}

Now, to the $_2\phi_1$'s in \eqref{2psi2gl} we apply Heine's~\cite{heine}
$q$-Euler transformation (see \cite[Eq.~(1.4.3)]{grhyp}):
\begin{equation}\label{qeulertf}
{}_2\phi_1\!\left[\begin{matrix}a,b\\c\end{matrix}\,;q,z\right]=
\frac{(abz/c;q)_{\infty}}{(z;q)_{\infty}}\,
{}_2\phi_1\!\left[\begin{matrix}c/a,c/b\\c\end{matrix}\,;q,
\frac {abz}c\right],
\end{equation}
where $\max(|z|,|abz/c|)<1$. Hence, by \eqref{qeulertf}
the expression in \eqref{2psi2gl} is transformed into
\begin{equation}\label{2psi2gl1}
\frac{(b_2/a_1,a_2,b_1/a_1,q,a_1z,q/a_1z;q)_{\infty}}
{(a_2/a_1,b_1,b_2,q/a_1,z,q/z;q)_{\infty}}
{}_2\phi_1\!\left[\begin{matrix}a_1q/b_1,a_1q/b_2\\
a_1q/a_2\end{matrix}\,;q,\frac{b_1b_2}{a_1a_2z}\right]
+\operatorname{idem}(a_1;a_2),
\end{equation}
which is exactly the right side of the $r=2$ case of \eqref{slgentfu}.

Now suppose that the formula in \eqref{slgentfu}
is already shown for any integer $t$ where $1\le t<r$.
To show the transformation for $t=r$,
we multiply both sides of \eqref{key32gl} by
\begin{equation*}
\frac{(a_3,\dots,a_r;q)_n}{(b_3,\dots,b_r;q)_n}z^n
\end{equation*}
and sum over all integers $n$.

On the right side we obtain
\begin{equation*}
{}_r\psi_r\!\left[\begin{matrix}a_1,a_2,\dots,a_r\\b_1,b_2,
\dots,b_r\end{matrix}\,;q,z\right].
\end{equation*}
On the left side we obtain
\begin{multline*}
\frac{(b_1b_2/a_2,b_2/a_1,a_2,b_1/a_1;q)_{\infty}}
{(a_2/a_1,b_1,b_1b_2/a_1a_2,b_2;q)_{\infty}}
\sum_{n=-\infty}^{\infty}\frac{(a_3,\dots,a_r;q)_n}{(b_3,\dots,b_r;q)_n}z^n\\
\times\sum_{k=0}^{\infty}\frac{(b_1/a_2,b_2/a_2;q)_k(a_1;q)_{n+k}}
{(q,a_1q/a_2;q)_k(b_1b_2/a_2;q)_{n+k}}q^k
+\operatorname{idem}(a_1;a_2)\\
=\frac{(b_1b_2/a_2,b_2/a_1,a_2,b_1/a_1;q)_{\infty}}
{(a_2/a_1,b_1,b_1b_2/a_1a_2,b_2;q)_{\infty}}
\sum_{k=0}^{\infty}\frac{(b_1/a_2,b_2/a_2,a_1;q)_k}
{(q,a_1q/a_2,b_1b_2/a_2;q)_k}q^k\\\times
\sum_{n=-\infty}^{\infty}\frac{(a_1q^k,a_3,\dots,a_r;q)_n}
{(b_1b_2q^k/a_2,b_3,\dots,b_r;q)_n}z^n
+\operatorname{idem}(a_1;a_2).
\end{multline*}
Now to the inner sums we apply the inductive hypothesis (i.e., the
$r\mapsto r-1$ case of \eqref{slgentfu}), and we obtain for the last
expression
\begin{multline*}
\frac{(b_1b_2/a_2,b_2/a_1,a_2,b_1/a_1;q)_{\infty}}
{(a_2/a_1,b_1,b_1b_2/a_1a_2,b_2;q)_{\infty}}
\sum_{k=0}^{\infty}\frac{(b_1/a_2,b_2/a_2,a_1;q)_k}
{(q,a_1q/a_2,b_1b_2/a_2;q)_k}q^k\\\times
\Bigg[\frac{(q,a_3,\dots,a_r,b_1b_2/a_1a_2,b_3q^{-k}/a_1,\dots,
b_rq^{-k}/a_1,a_1zq^k,q^{1-k}/a_1z;q)_{\infty}}
{(q^{1-k}/a_1,a_3q^{-k}/a_1,\dots,a_rq^{-k}/a_1,
b_1b_2q^k/a_2,b_3,\dots,b_r,z,q/z;q)_{\infty}}\\\times
\sum_{n=0}^{\infty}\frac{(a_1a_2q/b_1b_2,a_1q^{1+k}/b_3,\dots,
a_1q^{1+k}/b_r;q)_n}{(q,a_1q^{1+k}/a_3,\dots,a_1q^{1+k}/a_r;q)_n}
\left(\frac{b_1\dots b_r}{a_1\dots a_rz}\right)^n\\
+\Bigg(\frac{(q,a_1q^k,a_4,\dots,a_r,b_1b_2q^k/a_2a_3,b_3/a_3,\dots,
b_r/a_3,a_3z,q/a_3z;q)_{\infty}}
{(q/a_3,a_1q^k/a_3,a_4/a_3,\dots,a_r/a_3,
b_1b_2q^k/a_2,b_3,\dots,b_r,z,q/z;q)_{\infty}}\\\times
\sum_{n=0}^{\infty}\frac{(a_2a_3q^{1-k}/b_1b_2,a_3q/b_3,\dots,
a_3q/b_r;q)_n}{(q,a_3q^{1-k}/a_1,a_3q/a_4,\dots,a_3q/a_r;q)_n}
\left(\frac{b_1\dots b_r}{a_1\dots a_rz}\right)^n\\
+\operatorname{idem}(a_3;a_4,\dots,a_r)\Bigg)
\Bigg]+\operatorname{idem}(a_1;a_2),
\end{multline*}
which can be simplified to
\begin{multline}\label{slgen1}
\Bigg[\frac{(b_2/a_1,a_2,b_1/a_1,
q,a_3,\dots,a_r,b_3/a_1,\dots,
b_r/a_1,a_1z,q/a_1z;q)_{\infty}}
{(a_2/a_1,b_1,b_2,q/a_1,a_3/a_1,\dots,a_r/a_1,
b_3,\dots,b_r,z,q/z;q)_{\infty}}\\\times
\sum_{k=0}^{\infty}\frac{(b_1/a_2,b_2/a_2,a_1q/b_3,\dots,a_1q/b_r;q)_k}
{(q,a_1q/a_2,a_1q/a_3,\dots,a_1q/a_r;q)_k}
\left(\frac{b_3\dots b_rq}{a_3\dots a_rz}\right)^k\\\times
\sum_{n=0}^{\infty}\frac{(a_1a_2q/b_1b_2,a_1q^{1+k}/b_3,\dots,
a_1q^{1+k}/b_r;q)_n}{(q,a_1q^{1+k}/a_3,\dots,a_1q^{1+k}/a_r;q)_n}
\left(\frac{b_1\dots b_r}{a_1\dots a_rz}\right)^n\\
+\Bigg(\frac{(b_2/a_1,a_2,b_1/a_1,
q,a_1,a_4,\dots,a_r,b_1b_2/a_2a_3,b_3/a_3,\dots,
b_r/a_3,a_3z,q/a_3z;q)_{\infty}}
{(a_2/a_1,b_1,b_1b_2/a_1a_2,b_2,
q/a_3,a_1/a_3,a_4/a_3,\dots,a_r/a_3,
b_3,\dots,b_r,z,q/z;q)_{\infty}}\\\times
\sum_{k=0}^{\infty}\frac{(b_1/a_2,b_2/a_2,a_1/a_3;q)_k}
{(q,a_1q/a_2,b_1b_2/a_2a_3;q)_k}q^k\\\times
\sum_{n=0}^{\infty}\frac{(a_2a_3q^{1-k}/b_1b_2,a_3q/b_3,\dots,
a_3q/b_r;q)_n}{(q,a_3q^{1-k}/a_1,a_3q/a_4,\dots,a_3q/a_r;q)_n}
\left(\frac{b_1\dots b_r}{a_1\dots a_rz}\right)^n\\
+\operatorname{idem}(a_3;a_4,\dots,a_r)\Bigg)
\Bigg]+\operatorname{idem}(a_1;a_2).
\end{multline}

We have in \eqref{slgen1} a sum of $2(r-1)$ double sums.
Accordingly, for more clarity, let us write the whole expression
in \eqref{slgen1} as
\begin{align}\label{tu1}\notag
T_1(a_1,a_2)+{} &T_2(a_1,a_2)+\dots +T_{r-1}(a_1,a_2)\\
{}+ { } U_1(a_1,a_2)+{} &U_2(a_1,a_2)+\dots +U_{r-1}(a_1,a_2),
\end{align}
where, by definition of ``idem", $U_i(a_1,a_2)=T_i(a_2,a_1)$
for $i=1,\dots,r-1$. Further, $\sum_{i=2}^{r-1}T_i(a_1,a_2)
=T_2(a_1,a_2)+\operatorname{idem}(a_3;a_4,\dots,a_r)$
(and $\sum_{i=2}^{r-1}U_i(a_1,a_2)
=U_2(a_1,a_2)+\operatorname{idem}(a_3;a_4,\dots,a_r)$).
Below, we will selectively perform manipulations with the respective terms
$T_i$ and $U_i$ ($i=1,\dots,r-1$). 

To evaluate $T_1(a_1,a_2)$ (and hence also $U_1(a_1,a_2)$),
we first shift the index $n$ of the inner sum in $T_1(a_1,a_2)$ (or,
equivalently, in the first term in \eqref{slgen1}) by $-k$
and then interchange the double sum.
Symbolically, we apply the interchange of summations as in \eqref{intchs}.
Thus, we obtain (using some elementary identities for $q$-shifted factorials)
\begin{multline*}
T_1(a_1,a_2)=
\frac{(q,a_2,\dots,a_r,b_1/a_1,\dots,b_r/a_1,a_1z,q/a_1z;q)_{\infty}}
{(q/a_1,a_2/a_1,\dots,a_r/a_1,
b_1,\dots,b_r,z,q/z;q)_{\infty}}\\\times
\sum_{n=0}^{\infty}\frac{(a_1a_2q/b_1b_2,a_1q/b_3,\dots,
a_1q/b_r;q)_n}{(q,a_1q/a_3,\dots,a_1q/a_r;q)_n}
\left(\frac{b_1\dots b_r}{a_1\dots a_rz}\right)^n\\\times
\sum_{k=0}^n\frac{(b_1/a_2,b_2/a_2,q^{-n};q)_k}
{(q,a_1q/a_2,b_1b_2q^{-n}/a_1a_2;q)_k}q^k.
\end{multline*}
Now the inner sum can be evaluated by the terminating $q$-Pfaff--Saalsch\"utz
summation (cf.~\cite[Eq.~(II.12)]{grhyp}),
\begin{equation}\label{32gl}
{}_3\phi_2\!\left[\begin{matrix}
a,b,q^{-n}\\c,abq^{1-n}/c\end{matrix};\,q,q\right]
=\frac{(c/a,c/b;q)_n}{(c,c/ab;q)_n},
\end{equation}
which simplifies the last expression,
$T_1(a_1,a_2)$, to
\begin{multline}\label{slgen2}
\frac{(q,a_2,\dots,a_r,b_1/a_1,\dots,b_r/a_1,a_1z,q/a_1z;q)_{\infty}}
{(q/a_1,a_2/a_1,\dots,a_r/a_1,
b_1,\dots,b_r,z,q/z;q)_{\infty}}\\\times
\sum_{n=0}^{\infty}\frac{(a_1q/b_1,a_1q/b_2,\dots,
a_1q/b_r;q)_n}{(q,a_1q/a_2,\dots,a_1q/a_r;q)_n}
\left(\frac{b_1\dots b_r}{a_1\dots a_rz}\right)^n.
\end{multline}

Next, we consider $T_2(a_1,a_2)$ and $U_2(a_1,a_2)$.
By interchanging the double sums in $T_2(a_1,a_2)$ and in $U_2(a_1,a_2)$
we obtain
\begin{multline}\label{slgen3}
T_2(a_1,a_2)+U_2(a_1,a_2)\\
=\frac{(b_2/a_1,a_2,b_1/a_1,
q,a_1,a_4,\dots,a_r,b_1b_2/a_2a_3,b_3/a_3,\dots,
b_r/a_3,a_3z,q/a_3z;q)_{\infty}}
{(a_2/a_1,b_1,b_1b_2/a_1a_2,b_2,
q/a_3,a_1/a_3,a_4/a_3,\dots,a_r/a_3,
b_3,\dots,b_r,z,q/z;q)_{\infty}}\\\times
\sum_{n=0}^{\infty}\frac{(a_2a_3q/b_1b_2,a_3q/b_3,\dots,
a_3q/b_r;q)_n}{(q,a_3q/a_1,a_3q/a_4,\dots,a_3q/a_r;q)_n}
\left(\frac{b_1\dots b_r}{a_1\dots a_rz}\right)^n\\\times
\sum_{k=0}^{\infty}\frac{(b_1/a_2,b_2/a_2,a_1q^{-n}/a_3;q)_k}
{(q,a_1q/a_2,b_1b_2q^{-n}/a_2a_3;q)_k}q^k
+\operatorname{idem}(a_1;a_2).
\end{multline}
Now, to the inner sum of the {\em first} double sum (but not of the second!)
in \eqref{slgen3} we apply the nonterminating $_3\phi_2$ summation in
\eqref{32ntgl}, i.e.,
specifically we apply
\begin{multline}\label{32ntgl1}
\sum_{k=0}^{\infty}\frac{(b_1/a_2,b_2/a_2,a_1q^{-n}/a_3;q)_k}
{(q,a_1q/a_2,b_1b_2q^{-n}/a_2a_3;q)_k}q^k\\
=\frac{(a_2/a_1,b_2q^{-n}/a_3,b_1q^{-n}/a_3,b_1b_2/a_1a_2;q)_{\infty}}
{(b_1/a_1,b_2/a_2,a_2q^{-n}/a_3,b_1b_2q^{-n}/a_2a_3;q)_{\infty}}\\
-\frac{(a_2/a_1,b_1/a_2,b_2/a_2,a_1q^{-n}/a_3,
b_1b_2q^{-n}/a_1a_3;q)_{\infty}}
{(a_1/a_2,b_1/a_1,b_2/a_1,a_2q^{-n}/a_3,
b_1b_2q^{-n}/a_2a_3;q)_{\infty}}\\\times
\sum_{k=0}^{\infty}\frac{(b_1/a_1,b_2/a_1,a_2q^{-n}/a_3;q)_k}
{(q,a_2q/a_1,b_1b_2q^{-n}/a_1a_3;q)_k}q^k.
\end{multline}
The result of the application of this (two term) summation is that
the first term in \eqref{slgen3}, $T_2(a_1,a_2)$, is split into
two parts, one single sum $T_2^{\prime}(a_1,a_2)$ and one double
sum $T_2^{\prime\prime}(a_1,a_2)$.
Formally, we have
$$
T_2(a_1,a_2)+U_2(a_1,a_2)=[T_2^{\prime}(a_1,a_2)+
T_2^{\prime\prime}(a_1,a_2)]+U_2(a_1,a_2).
$$
But $T_2^{\prime\prime}(a_1,a_2)$ is precisely $-U_2(a_1,a_2)$
(as can be readily checked), so two terms cancel, thus
$$
T_2(a_1,a_2)+U_2(a_1,a_2)=T_2^{\prime}(a_1,a_2).
$$
Now, the last expression, $T_2^{\prime}(a_1,a_2)$, can be simplified to
\begin{multline}\label{slgen4}
\frac{(q,a_1,a_2,a_4,\dots,a_r,b_1/a_3,\dots,
b_r/a_3,a_3z,q/a_3z;q)_{\infty}}
{(q/a_3,a_1/a_3,a_2/a_3,a_4/a_3,\dots,a_r/a_3,
b_1,\dots,b_r,z,q/z;q)_{\infty}}\\\times
\sum_{n=0}^{\infty}\frac{(a_3q/b_1,a_3q/b_2,\dots,
a_3q/b_r;q)_n}{(q,a_3q/a_1,a_3q/a_2,a_3q/a_4,\dots,a_3q/a_r;q)_n}
\left(\frac{b_1\dots b_r}{a_1\dots a_rz}\right)^n.
\end{multline}
It is easy to see that \eqref{slgen4} equals
\eqref{slgen2} where $a_1$ and $a_3$ are interchanged.
Collecting all terms, according to \eqref{tu1}, completes our derivation
of Slater's $_r\psi_r$ transformation formula
in \eqref{slgentfu}.

\medskip
We conclude this section considering Slater's $_r\psi_r$
transformation in its general form:
\begin{multline}\label{slgentf}
{}_r\psi_r\!\left[\begin{matrix}a_1,a_2,\dots,a_r\\b_1,b_2,
\dots,b_r\end{matrix}\,;q,z\right]\\
=\frac{(\frac{c_1}{a_1},\dots,\frac{c_1}{a_r},c_2,\dots,c_r,
\frac q{c_2},\dots,\frac q{c_r},\frac{b_1q}{c_1},\dots,\frac{b_rq}{c_1},
Ac_1z,\frac q{Ac_1z};q)_{\infty}}
{(\frac q{a_1},\dots,\frac q{a_r},\frac{c_1}{c_2},\dots,\frac{c_1}{c_r},
\frac{c_2q}{c_1},\dots,\frac{c_rq}{c_1},
b_1,\dots,b_r,Azq,\frac 1{Az};q)_{\infty}}\\\times
{}_r\psi_r\!\left[\begin{matrix}\frac {a_1q}{c_1},
\dots,\frac {a_rq}{c_1}\\
\frac {b_1q}{c_1},\dots,\frac{b_rq}{c_1}\end{matrix}\,;q,z\right]
+\operatorname{idem}(c_1;c_2,\dots,c_r),
\end{multline}
where $A=a_1\dots a_r/c_1\dots c_r$, and
the series either terminate, or $|b_1\dots b_r/a_1\dots a_r|<|z|<1$,
for convergence.

The transformation in \eqref{slgentf}
which involves only bilateral series is much more general than
the transformation in \eqref{slgentfu} since in \eqref{slgentf}
we have $r$ additional parameters $c_1,\dots,c_r$. 
It is not difficult to see that the special case $c_i=a_iq$, $i=1,\dots,r$,
of \eqref{slgentf} is exactly \eqref{slgentfu}.
However, for some purposes we rather
consider a transformation equivalent to \eqref{slgentf},
see \eqref{slgentf1}, below.

The natural question arises whether we can also prove the more general
transformation \eqref{slgentf}, which involves only bilateral series,
by elementary means.
Unfortunately, we were not able to derive \eqref{slgentf} directly by
the method of this article. Instead, we can extend \eqref{slgentfu}
to \eqref{slgentf} by an $r$-fold application of Ismail's argument.
This works similar as in Section~\ref{sec2r} where we described
how the transformation in \eqref{slw2r2rgl} for well-poised
$_{2r}\psi_{2r}$ series can be extended to the transformation in
\eqref{sl2r2rgl}.
We provide a sketch of how Ismail's argument is applied here.

First, let us transform the identity in \eqref{slgentf} to an equivalent
identity by replacing $z$ by $z/A$, shifting the summation indices
by one, and reversing the infinite sums on the right hand side.
Specifically, we apply
\begin{multline*}
\sum_{k=-\infty}^{\infty}\frac{(\frac{a_1q}{c_1},\dots,\frac{a_rq}{c_1};q)_k}
{(\frac{b_1q}{c_1},\dots,\frac{b_rq}{c_1};q)_k}
\left(\frac{c_1\dots c_rz}{a_1\dots a_r}\right)^k\\
=\sum_{k=-\infty}^{\infty}\frac{(\frac{a_1q}{c_1},\dots,
\frac{a_rq}{c_1};q)_{-1-k}}
{(\frac{b_1q}{c_1},\dots,\frac{b_rq}{c_1};q)_{-1-k}}
\left(\frac{c_1\dots c_rz}{a_1\dots a_r}\right)^{-1-k}\\
=\frac{(1-\frac{c_1}{b_1})\dots(1-\frac{c_1}{b_r})}
{(1-\frac{c_1}{a_1})\dots(1-\frac{c_1}{a_r})}
\frac{b_1\dots b_r}{c_1\dots c_rz}
\sum_{k=-\infty}^{\infty}\frac{(\frac{c_1q}{b_1},\dots,\frac{c_1q}{b_r};q)_k}
{(\frac{c_1q}{a_1},\dots,\frac{c_1q}{a_r};q)_k}
\left(\frac{b_1\dots b_r}{a_1\dots a_rz}\right)^k
\end{multline*}
to the $_r\psi_r$'s on the right hand side of \eqref{slgentf}.
Hence, \eqref{slgentf} becomes
\begin{multline}\label{slgentf1}
{}_r\psi_r\!\left[\begin{matrix}a_1,a_2,\dots,a_r\\b_1,b_2,
\dots,b_r\end{matrix}\,;q,\frac{c_1\dots c_rz}{a_1\dots a_r}\right]\\
=\frac{(\frac{c_1q}{a_1},\dots,\frac{c_1q}{a_r},c_2,\dots,c_r,
\frac q{c_2},\dots,\frac q{c_r},\frac{b_1}{c_1},\dots,\frac{b_r}{c_1},
c_1z,\frac q{c_1z};q)_{\infty}}
{(\frac q{a_1},\dots,\frac q{a_r},\frac{c_1q}{c_2},\dots,\frac{c_1q}{c_r},
\frac{c_2}{c_1},\dots,\frac{c_r}{c_1},
b_1,\dots,b_r,z,\frac qz;q)_{\infty}}\\\times
{}_r\psi_r\!\left[\begin{matrix}\frac {c_1q}{b_1},
\dots,\frac{c_1q}{b_r}\\
\frac {c_1q}{a_1},\dots,\frac{c_1q}{a_r}\end{matrix}\,;q,
\frac{b_1\dots b_r}{c_1\dots c_rz}\right]
+\operatorname{idem}(c_1;c_2,\dots,c_r),
\end{multline}
where the series either terminate, or $|b_1\dots b_r/a_1\dots a_r|<|z|<1$,
for convergence.

In \eqref{slgentf1}, Slater's general $_r\psi_r$ transformation is written 
in a more convenient form for us since here we immediately see that
the special case $c_i=a_i$, $i=1,\dots,r$, of \eqref{slgentf1}
is exactly \eqref{slgentfu}. In the following, we closely follow the
final analysis of Section~\ref{sec2r}.

It is easy to see that both sides of \eqref{slgentf1} are analytic in
each $a_1^{-1},a_2^{-1},\dots,a_r^{-1}$ in a domain around the origin.
We know that the identity is true when $a_i=c_i$, for $i=1,\dots,r$.
What needs to be done is to extend \eqref{slgentfu} first by an additional
parameter $a_1$, then by $a_2$, etc. This means that if we have
already extended \eqref{slgentfu} by $a_1,\dots,a_j$, what we
should have derived is the identity in \eqref{slgentf1} where $a_i=c_i$,
for $i=j+1,\dots,r$. So, we proceed by induction starting with $j=0$
(identity \eqref{slgentfu}) and ending with $j=r$ (identity
\eqref{slgentf1}). In the inductive step, we consider the
transformation \eqref{slgentf1} where $a_i=c_i$, for $i=j+1,\dots,r$.
We call that identity $E_{j+1}$. We need to show that $E_{j+1}$ is true,
provided $E_j$ is true, which is \eqref{slgentf1} where $a_i=c_i$,
for $i=j,\dots,r$.
We observe that both sides of $E_{j+1}$ are
analytic in $1/a_j$ around the origin.
The identity is true for $a_j=c_jq^{-n_j}$, $n_j=0,1,2,\dots$.
This follows by the $c_j\mapsto c_jq^{-n_j}$ case
of the inductive hypothesis, $E_j$. (This can be verified by looking at all
the terms involving $j$. Further note that the index of the
$j$-th sum is shifted by $-n_j$, since the $j$-th bilateral series
becomes a unilateral series.) Since $E_{j+1}$ is true for
an infinite sequence of $1/a_j$ which has an accumulation point,
namely 0, in the interior of the domain $\mathcal D$ of analyticity of
$1/a_j$, we can apply the identity theorem to deduce that $E_{j+1}$
is true for $1/a_j$ throughout $\mathcal D$.
Now, by induction, \eqref{slgentf1} follows, with the general additional
parameters $a_1,\dots,a_r$.

\section{Transformations of Chu--Gasper--Karlsson--Minton-type}\label{secc}

We say that a basic hypergeometric series is of
{\em Gasper--Karlsson--Minton-type} if there are $s$ upper parameters
$a_1,\dots,a_s$ and $s$ lower parameters $b_1,\dots,b_s$ such that
each $a_i$ differs from $b_i$ multiplicatively
by a nonnegative integer power of $q$,
i.e\@. $a_i=b_iq^{m_i}$, $m_i\ge 0$ (see~\cite[Sec.~1.9]{grhyp}).
Originally, Minton~\cite{minton} and Karlsson~\cite{karlsson} had
discovered some corresponding summation formulae for ordinary
hypergeometric series (where there are $s$ upper parameters
$a_1,\dots,a_s$ and $s$ lower parameters $b_1,\dots,b_s$ such that
each $a_i$ differs from $b_i$ additively by a nonnegative integer,
i.e\@. $a_i=b_i+m_i$, $m_i\ge 0$).
Gasper~\cite{gassum} found $q$-analogues of
Karlsson and Minton's summations and also extended these to
transformation formulae. Later, Chu~\cite{chubs},\cite{chuwp} found
bilateral summations and transformations of Gasper--Karlsson--Minton-type
which included all the earlier known identities of
Gasper--Karlsson--Minton-type as special cases. Accordingly,
we call such bilateral identities to be of
{\em Chu--Gasper--Karlsson--Minton-type}.

Using here, and in the following, the notation $\{x_{\nu}\}$ for the
$s$ basic hypergeometric parameters $x_1,\dots,x_s$ ($\nu=1,\dots,s$),
and also $|m|=\sum_{i=1}^s m_i$, for brevity,
Chu's~\cite[Eq.~(15)]{chubs} transformation formula
for a specific $_{2+s}\psi_{2+s}$ series of
Chu--Gasper--Karlsson--Minton-type can be stated as follows:
\begin{multline}\label{chutf}
{}_{2+s}\psi_{2+s}\!\left[\begin{matrix}a,b,\{h_{\nu}q^{m_{\nu}}\}\\
c,d,\{h_{\nu}\}\end{matrix}\,;q,
\frac{q^{1-N}}a\right]
=b^N\frac{(q,bq/a,c/b,d/b;q)_{\infty}}
{(q/a,q/b,c,d;q)_{\infty}}
\prod_{i=1}^s\frac{(h_i/b;q)_{m_i}}{(h_i;q)_{m_i}}\\\times
{}_{2+s}\phi_{1+s}\!\left[\begin{matrix}bq/c,bq/d,\{bq/h_{\nu}\}\\
bq/a,\{bq^{1-m_{\nu}}/h_{\nu}\}\end{matrix}\,;q,
\frac{cdq^{N-|m|-1}}b\right],
\end{multline}
where $N$ is an arbitrary integer,
and where the series either terminate, or $|q/a|<|q^N|<|bq^{|m|+1}/cd|$,
for convergence.

On the other hand, Chu's~\cite[Theorem~2]{chuwp} summation
formula for a specific very-well-poised $_{6+2s}\psi_{6+2s}$ series of
Chu--Gasper--Karlsson--Minton-type is
\begin{multline}\label{chugl}
{}_{6+2s}\psi_{6+2s}\!\left[\begin{matrix}q\sqrt{a},-q\sqrt{a},b,c,d,\frac ad,
\{h_{\nu}\},\{\frac{aq^{1+m_{\nu}}}{h_{\nu}}\}\\
\sqrt{a},-\sqrt{a},\frac{aq}b,\frac{aq}c,\frac{aq}d,dq,
\{\frac{aq}{h_{\nu}}\},\{h_{\nu}q^{-m_{\nu}}\}\end{matrix}\,;q,
\frac{aq^{1-|m|}}{bc}\right]\\
=\frac{(q,q,aq,q/a,aq/bd,aq/cd,dq/b,dq/c;q)_{\infty}}
{(aq/b,aq/c,aq/d,dq/a,q/b,q/c,q/d,dq;q)_{\infty}}
\prod_{i=1}^s\frac{(aq/dh_i,dq/h_i;q)_{m_i}}{(aq/h_i,q/h_i;q)_{m_i}},
\end{multline}
provided the series either terminates, or $|aq^{1-|m|}/bc|<1$,
for convergence.

Chu established both of these formulae
by means of partial fraction expansions.
Haglund~\cite[pp.~415--416]{hagl} noticed that the transformation
in \eqref{chutf}
can be obtained by specializing Slater's general transformation
\eqref{slgentf} for $_r\psi_r$ series. It is also true that the
summation in \eqref{chugl} can be derived by specializing Slater's
transformation \eqref{sl2r2rgl}
for well-poised $_{2r}\psi_{2r}$ series.

Rather then specializing Slater's transformations,
in the spirit of this article we give elementary derivations of
Chu--Gasper--Karlsson--Minton-type identities.
We show here that both of Slater's general transformation
formulae, \eqref{slgentf} and \eqref{sl2r2rgl}, can be extended by induction
to transformations of Chu--Gasper--Karlsson--Minton-type,
see Propositions~\ref{chugen1} and \ref{chugen2} below.
The simple analysis involves interchanging of sums,
and the expansion of certain factors in terms
of the $q$-binomial theorem~\eqref{10tgl} or $q$-Pfaff--Saalsch\"utz
summation~\eqref{32gl}, respectively.
Finally, we state some interesting special cases
(in addition to \eqref{chutf} and \eqref{chugl}) of the general
Chu--Gasper--Karlsson--Minton-type transformations explicitly, see
Corollaries~\ref{km1}, \ref{km2} and \ref{km3}.

\begin{Proposition}[A general $_{r+s}\psi_{r+s}$
Chu--Gasper--Karlsson--Minton-type transformation]\label{chugen1}
Let $a_1,\dots,a_r$, $b_1,\dots,b_r$, $c_1,\dots,c_r$,
$h_1,\dots,h_s$, and $z$ be indeterminate, let $m_1,\dots,m_s$
be nonnegative integers, let $A=a_1\dots a_r/c_1\dots c_r$,
$|m|=\sum_{i=1}^s m_i$, and suppose
that the series in \eqref{chugen1gl} are well-defined. Then
\begin{multline}\label{chugen1gl}
{}_{r+s}\psi_{r+s}\!\left[\begin{matrix}a_1,\dots,a_r,\{h_{\nu}q^{m_{\nu}}\}\\
b_1,\dots,b_r,\{h_{\nu}\}\end{matrix}\,;q,z\right]
=\prod_{i=1}^s\frac{(\frac{h_iq}{c_1};q)_{m_i}}{(h_i;q)_{m_i}}\\\times
\frac{(\frac{c_1}{a_1},\dots,\frac{c_1}{a_r},c_2,\dots,c_r,
\frac q{c_2},\dots,\frac q{c_r},\frac{b_1q}{c_1},\dots,\frac{b_rq}{c_1},
Ac_1z,\frac q{Ac_1z};q)_{\infty}}
{(\frac q{a_1},\dots,\frac q{a_r},\frac{c_1}{c_2},\dots,\frac{c_1}{c_r},
\frac{c_2q}{c_1},\dots,\frac{c_rq}{c_1},
b_1,\dots,b_r,Azq,\frac 1{Az};q)_{\infty}}\\\times
{}_r\psi_r\!\left[\begin{matrix}\frac {a_1q}{c_1},
\dots,\frac {a_rq}{c_1},\{\frac{h_{\nu}q^{1+m_{\nu}}}{c_1}\}\\
\frac {b_1q}{c_1},\dots,\frac{b_rq}{c_1},
\{\frac{h_{\nu}q}{c_1}\}\end{matrix}\,;q,z\right]
+\operatorname{idem}(c_1;c_2,\dots,c_r),
\end{multline}
where the series either terminate,
or $|b_1\dots b_rq^{-|m|}/a_1\dots a_r|<|z|<1$.
\end{Proposition}

\begin{proof}
We proceed by induction on $s$. For $s=0$ the transformation is true
by Slater's general transformation~\eqref{slgentf}. So, suppose
the identity is already shown for $s\mapsto s-1$.
Then,
\begin{multline}\label{eq1}
{}_{r+s}\psi_{r+s}\!\left[\begin{matrix}a_1,\dots,a_r,\{h_{\nu}q^{m_{\nu}}\}\\
b_1,\dots,b_r,\{h_{\nu}\}\end{matrix}\,;q,z\right]\\
=\sum_{k=-\infty}^{\infty}\frac{(a_1,\dots,a_r;q)_k}{(b_1,\dots,b_r;q)_k}
\prod_{i=1}^{s-1}\frac{(h_iq^{m_i};q)_k}{(h_i;q)_k}z^k\,
\frac{(h_sq^{m_s};q)_k}{(h_s;q)_k}\\
=\frac 1{(h_s;q)_{m_s}}
\sum_{k=-\infty}^{\infty}\frac{(a_1,\dots,a_r;q)_k}{(b_1,\dots,b_r;q)_k}
\prod_{i=1}^{s-1}\frac{(h_iq^{m_i};q)_k}{(h_i;q)_k}z^k\,(h_sq^k;q)_{m_s}.
\end{multline}
Now we expand the last $q$-shifted factorial
by the terminating $q$-binomial theorem
(cf.~\cite[Eq.~(II.4)]{grhyp}),
\begin{equation}\label{10tgl}
{}_1\phi_0\!\left[\begin{matrix}q^{-n}\\
-\end{matrix}\,;q,z\right]=(zq^{-n};q)_n,
\end{equation}
which is just the $a\mapsto q^{-n}$ case of \eqref{10gl}.
That is, we apply
\begin{equation*}
(h_sq^k;q)_{m_s}=\sum_{j=0}^{m_s}\frac{(q^{-m_s};q)_j}{(q;q)_j}
\left(h_sq^{k+m_s}\right)^j
\end{equation*}
and obtain for the expression in \eqref{eq1}
\begin{multline}\label{eq2}
\frac 1{(h_s;q)_{m_s}}
\sum_{k=-\infty}^{\infty}\frac{(a_1,\dots,a_r;q)_k}{(b_1,\dots,b_r;q)_k}
\prod_{i=1}^{s-1}\frac{(h_iq^{m_i};q)_k}{(h_i;q)_k}z^k
\sum_{j=0}^{m_s}\frac{(q^{-m_s};q)_j}{(q;q)_j}
\left(h_sq^{k+m_s}\right)^j\\
=\frac 1{(h_s;q)_{m_s}}\sum_{j=0}^{m_s}\frac{(q^{-m_s};q)_j}{(q;q)_j}
\left(h_sq^{m_s}\right)^j\\\times
\sum_{k=-\infty}^{\infty}\frac{(a_1,\dots,a_r;q)_k}{(b_1,\dots,b_r;q)_k}
\prod_{i=1}^{s-1}\frac{(h_iq^{m_i};q)_k}{(h_i;q)_k}\left(zq^j\right)^k.
\end{multline}
Note that $j$ is bounded by $m_s$, thus the interchange of summations
in \eqref{eq2} is justified (provided
$|b_1\dots b_rq^{-|m|}/a_1\dots a_r|<|z|<1$).
Now we can apply the inductive hypothesis to the inner sum which gives us
\begin{multline*}
\frac 1{(h_s;q)_{m_s}}\sum_{j=0}^{m_s}\frac{(q^{-m_s};q)_j}{(q;q)_j}
\left(h_sq^{m_s}\right)^j
\prod_{i=1}^{s-1}\frac{(h_iq/c_1;q)_{m_i}}{(h_i;q)_{m_i}}\,
\frac{(c_1/a_1,\dots,c_1/a_r;q)_{\infty}}
{(q/a_1,\dots,q/a_r;q)_{\infty}}\\\times
\frac{(c_2,\dots,c_r,
q/c_2,\dots,q/c_r,b_1q/c_1,\dots,b_rq/c_1,Ac_1zq^j,q^{1-j}/Ac_1z;q)_{\infty}}
{(c_1/c_2,\dots,c_1/c_r,
c_2q/c_1,\dots,c_rq/c_1,
b_1,\dots,b_r,Azq^{1+j},q^{-j}/Az;q)_{\infty}}\\\times
\sum_{k=-\infty}^{\infty}\frac{(a_1q/c_1,\dots,a_rq/c_1;q)_k}
{(b_1q/c_1,\dots,b_rq/c_1;q)_k}
\prod_{i=1}^{s-1}\frac{(h_iq^{1+m_i}/c_1;q)_k}{(h_iq/c_1;q)_k}
\left(zq^j\right)^k\\
+\operatorname{idem}(c_1;c_2,\dots,c_r).
\end{multline*}
In this expression, we again interchange summations, and obtain
\begin{multline*}
\frac{(c_1/a_1,\dots,c_1/a_r,c_2,\dots,c_r,
q/c_2,\dots,q/c_r,b_1q/c_1,\dots,b_rq/c_1,Ac_1z,q/Ac_1z;q)_{\infty}}
{(q/a_1,\dots,q/a_r,c_1/c_2,\dots,c_1/c_r,
c_2q/c_1,\dots,c_rq/c_1,
b_1,\dots,b_r,Azq,1/Az;q)_{\infty}}\\\times
\frac 1{(h_s;q)_{m_s}}
\prod_{i=1}^{s-1}\frac{(h_iq/c_1;q)_{m_i}}{(h_i;q)_{m_i}}
\sum_{k=-\infty}^{\infty}\frac{(a_1q/c_1,\dots,a_rq/c_1;q)_k}
{(b_1q/c_1,\dots,b_rq/c_1;q)_k}
\prod_{i=1}^{s-1}\frac{(h_iq^{1+m_i}/c_1;q)_k}{(h_iq/c_1;q)_k}z^k\\\times
\sum_{j=0}^{m_s}\frac{(q^{-m_s};q)_j}{(q;q)_j}
\left(\frac{h_sq^{1+k+m_s}}{c_1}\right)^j
+\operatorname{idem}(c_1;c_2,\dots,c_r).
\end{multline*}
We simplify the inner sums, according to \eqref{10tgl},
\begin{equation*}
\sum_{j=0}^{m_s}\frac{(q^{-m_s};q)_j}{(q;q)_j}
\left(\frac{h_sq^{1+k+m_s}}{c_1}\right)^j=(h_sq^{1+k}/c_1;q)_{m_s},
\end{equation*}
and eventually deduce the proposition.
\end{proof}

If, in Proposition~\ref{chugen1}, we set $r=2$, $a_1\mapsto a$,
$a_2\mapsto b$, $b_1\mapsto c$, $b_2\mapsto d$, $c_1\mapsto e$,
and $z\mapsto eq^{-N}/ab$, then the second term on the right
side of \eqref{chugen1gl} vanishes. Furthermore, the
parameter $c_2$ cancels out in the first term on the right side.
We obtain
\begin{Corollary}\label{km1}
Let $a$, $b$, $c$, $d$, $e$, and $h_1,\dots,h_s$ be indeterminate,
let $N$ be an arbitrary integer, $m_1,\dots,m_s$
be nonnegative integers, and suppose that the series in
\eqref{chutfgen} are well-defined. Then
\begin{multline}\label{chutfgen}
{}_{2+s}\psi_{2+s}\!\left[\begin{matrix}a,b,\{h_{\nu}q^{m_{\nu}}\}\\
c,d,\{h_{\nu}\}\end{matrix}\,;q,\frac{eq^{-N}}{ab}\right]
=\left(\frac eq\right)^N\frac{(e/a,e/b,cq/e,dq/e;q)_{\infty}}
{(q/a,q/b,c,d;q)_{\infty}}\\\times
\prod_{i=1}^s\frac{(h_iq/e;q)_{m_i}}{(h_i;q)_{m_i}}\;
{}_{2+s}\psi_{2+s}\!\left[\begin{matrix}aq/e,bq/e,\{h_{\nu}q^{1+m_{\nu}}/e\}\\
cq/e,dq/e,\{h_{\nu}q/e\}\end{matrix}\,;q,\frac{eq^{-N}}{ab}\right],
\end{multline}
where the series either terminate, or $|e/ab|<|q^N|<|eq^{|m|}/cd|$,
for convergence.
\end{Corollary}

For $e=d$ the ${}_{2+s}\psi_{2+s}$ on the right side of \eqref{chutfgen}
reduces to a ${}_{2+s}\phi_{1+s}$ series. This gives a transformation
for a ${}_{2+s}\psi_{2+s}$ into a (multiple of a) ${}_{2+s}\phi_{1+s}$
series which is different from Chu's transformation in \eqref{chutf}.

We can also reverse the ${}_{2+s}\psi_{2+s}$ series on the right side of
\eqref{chutfgen}. We obtain
\begin{multline}\label{chutfgen1}
{}_{2+s}\psi_{2+s}\!\left[\begin{matrix}a,b,\{h_{\nu}q^{m_{\nu}}\}\\
c,d,\{h_{\nu}\}\end{matrix}\,;q,\frac{eq^{-N}}{ab}\right]
=\left(\frac eq\right)^N\frac{(e/a,e/b,cq/e,dq/e;q)_{\infty}}
{(q/a,q/b,c,d;q)_{\infty}}\\\times
\prod_{i=1}^s\frac{(h_iq/e;q)_{m_i}}{(h_i;q)_{m_i}}\;
{}_{2+s}\psi_{2+s}\!\left[\begin{matrix}e/c,e/d,\{e/h_{\nu}\}\\
e/a,e/b,\{eq^{-m_{\nu}}/h_{\nu}\}\end{matrix}\,;q,\frac{cdq^{N-|m|}}e\right],
\end{multline} 
where the series either terminate, or $|e/ab|<|q^N|<|eq^{|m|}/cd|$,
for convergence.
We immediately see that the special case $e=bq$ of \eqref{chutfgen1}
is exactly Chu's transformation~\eqref{chutf}.

Another noteworthy specialization of Proposition~\ref{chugen1}
is simply the $r=1$ case, rewritten in the following corollary:
\begin{Corollary}\label{km2}
Let $a$, $b$, $c$, $z$, and $h_1,\dots,h_s$ be indeterminate,
let $m_1,\dots,m_s$ be nonnegative integers, and suppose that the series in
\eqref{km2gl} are well-defined. Then
\begin{multline}\label{km2gl}
{}_{1+s}\psi_{1+s}\!\left[\begin{matrix}a,\{h_{\nu}q^{m_{\nu}}\}\\
b,\{h_{\nu}\}\end{matrix}\,;q,z\right]
=\frac{(c/a,bq/c,az,q/az;q)_{\infty}}
{(q/a,b,azq/c,c/az;q)_{\infty}}\\\times
\prod_{i=1}^s\frac{(h_iq/c;q)_{m_i}}{(h_i;q)_{m_i}}\;
{}_{1+s}\psi_{1+s}\!\left[\begin{matrix}aq/c,\{h_{\nu}q^{1+m_{\nu}}/c\}\\
bq/c,\{h_{\nu}q/c\}\end{matrix}\,;q,z\right],
\end{multline}
where the series either terminate, or $|bq^{-|m|}/a|<|z|<1$,
for convergence.
\end{Corollary}

We can set $c=b$ or $c=aq$ in \eqref{km2gl} to reduce the ${}_{1+s}\psi_{1+s}$
series on the right side to a ${}_{1+s}\phi_{s}$ series if we want.
If we first let $c\mapsto b$, and then $b\mapsto aq$ in \eqref{km2gl},
the series on the right side of \eqref{km2gl} reduces to 1 and we would
obtain the following summation:
\begin{equation}\label{km2rgl}
{}_{1+s}\psi_{1+s}\!\left[\begin{matrix}a,\{h_{\nu}q^{m_{\nu}}\}\\
aq,\{h_{\nu}\}\end{matrix}\,;q,z\right]
=\frac{(q,q,az,q/az;q)_{\infty}}
{(aq,q/a,z,q/z;q)_{\infty}}
\prod_{i=1}^s\frac{(h_iq/a;q)_{m_i}}{(h_i;q)_{m_i}},
\end{equation}
where the series either terminates, or $|q^{1-|m|}|<|z|<1$, for convergence.
Now, in the beginning of this article we assumed $|q|<1$. It follows from
the convergence condition ($|q^{1-|m|}|<1$) above that $|m|$ needs to be
either 0 (thus there are no parameters $h_i$
and we are left with a plain $_1\psi_1$), or that the series terminates,
and in our case, this would mean that the series reduces just to one term.
On the other hand, we may consider $|q|>1$ for 
a nonempty region of convergence.
In this case the convergence condition for
the bilateral series on the left side of \eqref{km2rgl} is
$1<|z|<|q^{1-|m|}|$. But unfortunately, our derivation of the identity
is not valid in this case. Clearly, if $|q|>1$, we could not have used
induction because the inductive basis would already have been false.
In total, \eqref{km2rgl} does not give a new summation. 

Nevertheless, we can still specialize Corollary~\ref{km2} to a summation
by choosing $c\mapsto b$, and then $b\mapsto aq^{1+n}$ where $n$
is a small positive integer. In this case, the series on the right side
of \eqref{km2gl} has only a finite number of terms and can be summed
explicitly. For instance, (when $s=n=m=1$) we have
\begin{equation}\label{ul1}
{}_2\psi_2\!\left[\begin{matrix}a,hq\\
aq^2,h\end{matrix}\,;q,z\right]
=\frac{(1-h/aq)-(1-h/a)z/q}{(1-h)}
\frac{(q^2,q,az,q/az;q)_{\infty}}
{(q/a,aq^2,z/q,q^2/z;q)_{\infty}},
\end{equation}
provided $|q|<|z|<1$. For $z=-h/a$ the right side of \eqref{ul1}
factors (completely into linear factors), and we have
\begin{equation}\label{ul2}
{}_2\psi_2\!\left[\begin{matrix}a,hq\\
aq^2,h\end{matrix}\,;q,-\frac ha\right]
=\frac{(1-h^2/a^2q)}{(1-h)}
\frac{(q^2,q,-h,-q/h;q)_{\infty}}
{(q/a,aq^2,-h/aq,-aq^2/h;q)_{\infty}},
\end{equation}
provided $|q|<|h/a|<1$.

Let us now consider transformations of Chu--Gasper--Karlsson--Minton-type 
for well-poised basic series.
We have
\begin{Proposition}[A well-poised $_{2r+2s}\psi_{2r+2s}$
Chu--Gasper--Karlsson--Minton type transformation]\label{chugen2}
Let $a$, $a_1,\dots,a_r$, $b_1,\dots,b_r$, and $h_1,\dots,h_s$
be indeterminate, let $m_1,\dots,m_s$
be nonnegative integers, $|m|=\sum_{i=1}^s m_i$, and suppose
that the series in \eqref{chugen2gl} are well-defined. Then
\begin{multline}\label{chugen2gl}
{}_{2r+2s}\psi_{2r+2s}\!\left[\begin{matrix}b_1,\dots,b_{2r},
\{h_{\nu}\},\{\frac{aq^{1+m_{\nu}}}{h_{\nu}}\}\\
\frac{a q}{b_1},\dots,\frac{a q}{b_{2r}},
\{\frac{aq}{h_{\nu}}\},\{h_{\nu}q^{-m_{\nu}}\}\end{matrix}\,;q,
-\frac{a^rq^{r-|m|}}{b_1\dots b_{2r}}\right]\\
=\frac{(a,\frac q{a},a_2,\dots,a_r,\frac q{a_2},\dots,\frac q{a_r},
\frac{a_2}{a},\dots,\frac{a_r}{a},
\frac{aq}{a_2},\dots,\frac{aq}{a_r},\frac{a_1q}{b_1},
\dots,\frac{a_1q}{b_{2r}};q)_{\infty}}
{(\frac q{b_1},\dots,\frac q{b_{2r}},
\frac{a q}{b_1},\dots,\frac{a q}{b_{2r}},
\frac{a_2}{a_1},\dots,\frac{a_r}{a_1},\frac{a_1q}{a_2},\dots,
\frac{a_1q}{a_r},\frac{a_1a_2}{a},\dots,
\frac{a_1a_r}{a};q)_{\infty}}\\\times
\frac{(\frac{a q}{a_1b_1},\dots,\frac{a q}{a_1b_{2r}},
\frac{a_1}{\sqrt{a}},-\frac{a_1}{\sqrt{a}},\frac{\sqrt{a}q}{a_1},
-\frac{\sqrt{a}q}{a_1};q)_{\infty}}
{(\frac{aq}{a_1a_2},\dots,\frac{aq}{a_1a_r},
\frac{a_1^2}{a},\frac{aq}{a_1^2},
\frac q{\sqrt{a}},-\frac q{\sqrt{a}},\sqrt{a},-\sqrt{a};q)_{\infty}}
\prod_{i=1}^s\frac{(\frac{a_1q}{h_i},\frac{aq}{a_1h_i};q)_{m_i}}
{(\frac{aq}{h_i},\frac q{h_i};q)_{m_i}}\\\times
{}_{2r+2s}\psi_{2r+2s}\!\left[\begin{matrix}\frac{a_1b_1}{a},
\dots,\frac{a_1b_{2r}}{a},
\{\frac{a_1h_{\nu}}a\},\{\frac{a_1q^{1+m_{\nu}}}{h_{\nu}}\}\\
\frac{a_1q}{b_1},\dots,\frac{a_1q}{b_{2r}},
\{\frac{a_1q}{h_{\nu}}\},\{\frac{a_1h_{\nu}q^{-m_{\nu}}}a\}\end{matrix}\,;q,
-\frac{a^rq^{r-|m|}}{b_1\dots b_{2r}}\right]\\
+\operatorname{idem}(a_1;a_2,\dots,a_r),
\end{multline}
where the series either terminate, or $|a^rq^{r-|m|}/b_1\dots b_{2r}|<1$,
for convergence.
\end{Proposition}

\begin{proof}
We proceed by induction on $s$. For $s=0$ the transformation is true
by Slater's transformation~\eqref{sl2r2rgl} for well-poised series.
So, suppose the identity is already shown for $s\mapsto s-1$.
Then,
\begin{multline}\label{eq3}
{}_{2r+2s}\psi_{2r+2s}\!\left[\begin{matrix}b_1,\dots,b_{2r},
\{h_{\nu}\},\{\frac{aq^{1+m_{\nu}}}{h_{\nu}}\}\\
\frac{a q}{b_1},\dots,\frac{a q}{b_{2r}},
\{\frac{aq}{h_{\nu}}\},\{h_{\nu}q^{-m_{\nu}}\}\end{matrix}\,;q,
-\frac{a^rq^{r-|m|}}{b_1\dots b_{2r}}\right]\\
=\sum_{k=-\infty}^{\infty}\frac{(b_1,\dots,b_{2r};q)_k}
{(aq/b_1,\dots,aq/b_{2r};q)_k}
\prod_{i=1}^{s-1}\frac{(h_i,aq^{1+m_i}/h_i;q)_k}
{(aq/h_i,h_iq^{-m_i};q)_k}
\left(-\frac{a^rq^{r-(m_1+\dots+m_s)}}{b_1\dots b_{2r}}\right)^k\\\times
\frac{(h_s,aq^{1+m_s}/h_s;q)_k}
{(aq/h_s,h_sq^{-m_s};q)_k}
=\frac{(b_{2r}q/h_s,aq/b_{2r}h_s;q)_{m_s}}{(aq/h_s,q/h_s;q)_{m_s}}\\\times
\sum_{k=-\infty}^{\infty}\frac{(b_1,\dots,b_{2r};q)_k}
{(aq/b_1,\dots,aq/b_{2r};q)_k}
\prod_{i=1}^{s-1}\frac{(h_i,aq^{1+m_i}/h_i;q)_k}
{(aq/h_i,h_iq^{-m_i};q)_k}
\left(-\frac{a^rq^{r-(m_1+\dots+m_{s-1})}}{b_1\dots b_{2r}}\right)^k\\\times
\frac{(aq^{1+k}/h_s,q^{1-k}/h_s;q)_{m_s}}
{(b_{2r}q/h_s,aq/b_{2r}h_s;q)_{m_s}}.
\end{multline}
Now we expand the last quotient of $q$-shifted factorials by the
terminating $q$-Pfaff--Saalsch\"utz summation~\eqref{32gl} and obtain
for the expression in \eqref{eq3}
\begin{multline}\label{eq4}
\frac{(b_{2r}q/h_s,aq/b_{2r}h_s;q)_{m_s}}{(aq/h_s,q/h_s;q)_{m_s}}
\sum_{k=-\infty}^{\infty}\frac{(b_1,\dots,b_{2r};q)_k}
{(aq/b_1,\dots,aq/b_{2r};q)_k}
\prod_{i=1}^{s-1}\frac{(h_i,aq^{1+m_i}/h_i;q)_k}
{(aq/h_i,h_iq^{-m_i};q)_k}\\\times
\left(-\frac{a^rq^{r-(m_1+\dots+m_{s-1})}}{b_1\dots b_{2r}}\right)^k
\sum_{j=0}^{m_s}\frac{(b_{2r}q^{-k}/a,b_{2r}q^k,q^{-m_s};q)_j}
{(q,b_{2r}q/h_s,b_{2r}h_sq^{-m_s}/a;q)_j}q^j\\
=\frac{(b_{2r}q/h_s,aq/b_{2r}h_s;q)_{m_s}}{(aq/h_s,q/h_s;q)_{m_s}}
\sum_{j=0}^{m_s}\frac{(b_{2r}/a,b_{2r},q^{-m_s};q)_j}
{(q,b_{2r}q/h_s,b_{2r}h_sq^{-m_s}/a;q)_j}q^j\\\times
\sum_{k=-\infty}^{\infty}\frac{(b_1,\dots,b_{2r-1},b_{2r}q^j;q)_k}
{(aq/b_1,\dots,aq/b_{2r-1},aq^{1-j}/b_{2r};q)_k}\\\times
\prod_{i=1}^{s-1}\frac{(h_i,aq^{1+m_i}/h_i;q)_k}
{(aq/h_i,h_iq^{-m_i};q)_k}
\left(-\frac{a^rq^{r-(m_1+\dots+m_{s-1})-j}}{b_1\dots b_{2r}}\right)^k.
\end{multline}
Note that $j$ is bounded by $m_s$, thus the interchange of summations
in \eqref{eq4} is justified (provided $|a^rq^{r-|m|}/b_1\dots b_{2r}|<1$).
Now we can apply the inductive hypothesis to the inner sum which gives us
\begin{multline*}
\frac{(b_{2r}q/h_s,aq/b_{2r}h_s;q)_{m_s}}{(aq/h_s,q/h_s;q)_{m_s}}
\sum_{j=0}^{m_s}\frac{(b_{2r}/a,b_{2r},q^{-m_s};q)_j}
{(q,b_{2r}q/h_s,b_{2r}h_sq^{-m_s}/a;q)_j}q^j\\\times
\frac{(a,q/a,a_2,\dots,a_r,q/a_2,\dots,q/a_r;q)_{\infty}}
{(q/b_1,\dots,q/b_{2r-1},q^{1-j}/b_{2r},aq/b_1,\dots,aq/b_{2r-1},
aq^{1-j}/b_{2r};q)_{\infty}}\\\times
\frac{(a_2/a,\dots,a_r/a,aq/a_2,\dots,aq/a_r,a_1q/b_1,\dots,a_1q/b_{2r-1},
a_1q^{1-j}/b_{2r};q)_{\infty}}
{(a_2/a_1,\dots,a_r/a_1,a_1q/a_2,\dots,a_1q/a_r,a_1a_2/a,
\dots,a_1a_r/a;q)_{\infty}}\\\times
\frac{(aq/a_1b_1,\dots,aq/a_1b_{2r-1},aq^{1-j}/a_1b_{2r},
a_1/\sqrt{a},-a_1/\sqrt{a},\sqrt{a}q/a_1,
-\sqrt{a}q/a_1;q)_{\infty}}
{(aq/a_1a_2,\dots,aq/a_1a_r,a_1^2/a,aq/a_1^2,q/\sqrt{a},-q/\sqrt{a},
\sqrt{a},-\sqrt{a};q)_{\infty}}\\\times
\prod_{i=1}^{s-1}\frac{(a_1q/h_i,aq/a_1h_i;q)_{m_i}}
{(aq/h_i,q/h_i;q)_{m_i}}
\sum_{k=-\infty}^{\infty}\frac{(a_1b_1/a,\dots,a_1b_{2r-1}/a,
a_1b_{2r}q^j/a;q)_k}
{(a_1q/b_1,\dots,a_1q/b_{2r-1},a_1q^{1-j}/b_{2r};q)_k}\\\times
\prod_{i=1}^{s-1}\frac{(a_1h_i/a,a_1q^{1+m_i}/h_i;q)_k}
{(a_1q/h_i,a_1h_iq^{-m_i}/a;q)_k}
\left(-\frac{a^rq^{r-(m_1+\dots+m_{s-1})-j}}{b_1\dots b_{2r}}\right)^k
+\operatorname{idem}(a_1;a_2,\dots,a_r).
\end{multline*}
In this expression, we again interchange summations, and obtain
\begin{multline*}
\frac{(b_{2r}q/h_s,aq/b_{2r}h_s;q)_{m_s}}{(aq/h_s,q/h_s;q)_{m_s}}
\frac{(a,q/a,a_2,\dots,a_r,q/a_2,\dots,q/a_r;q)_{\infty}}
{(q/b_1,\dots,q/b_{2r},aq/b_1,\dots,aq/b_{2r};q)_{\infty}}\\\times
\frac{(a_2/a,\dots,a_r/a,aq/a_2,\dots,aq/a_r,a_1q/b_1,\dots,
a_1q/b_{2r};q)_{\infty}}
{(a_2/a_1,\dots,a_r/a_1,a_1q/a_2,\dots,a_1q/a_r,a_1a_2/a,
\dots,a_1a_r/a;q)_{\infty}}\\\times
\frac{(aq/a_1b_1,\dots,aq/a_1b_{2r},a_1/\sqrt{a},-a_1/\sqrt{a},\sqrt{a}q/a_1,
-\sqrt{a}q/a_1;q)_{\infty}}
{(aq/a_1a_2,\dots,aq/a_1a_r,a_1^2/a,aq/a_1^2,q/\sqrt{a},-q/\sqrt{a},
\sqrt{a},-\sqrt{a};q)_{\infty}}\\\times
\prod_{i=1}^{s-1}\frac{(a_1q/h_i,aq/a_1h_i;q)_{m_i}}
{(aq/h_i,q/h_i;q)_{m_i}}\;
\sum_{k=-\infty}^{\infty}\frac{(a_1b_1/a,\dots,a_1b_{2r}/a;q)_k}
{(a_1q/b_1,\dots,a_1q/b_{2r};q)_k}\\\times
\prod_{i=1}^{s-1}\frac{(a_1h_i/a,a_1q^{1+m_i}/h_i;q)_k}
{(a_1q/h_i,a_1h_iq^{-m_i}/a;q)_k}
\left(-\frac{a^rq^{r-(m_1+\dots+m_{s-1})}}{b_1\dots b_{2r}}\right)^k\\\times
\sum_{j=0}^{m_s}\frac{(b_{2r}q^{-k}/a_1,a_1b_{2r}q^k/a,q^{-m_s};q)_j}
{(q,b_{2r}q/h_s,b_{2r}h_sq^{-m_s}/a;q)_j}q^j
+\operatorname{idem}(a_1;a_2,\dots,a_r).
\end{multline*}
We simplify the last inner sum, according to \eqref{32gl},
\begin{equation*}
\sum_{j=0}^{m_s}\frac{(b_{2r}q^{-k}/a_1,a_1b_{2r}q^k/a,
q^{-m_s};q)_j}{(q,b_{2r}q/h_s,b_{2r}h_sq^{-m_s}/a;q)_j}q^j
=\frac{(a_1q^{1+k}/h_s,aq^{1-k}/a_1h_s;q)_{m_s}}
{(b_{2r}q/h_s,aq/b_{2r}h_s;q)_{m_s}},
\end{equation*}
and eventually deduce the proposition.
\end{proof}

It is worth noting the case $a_1=b_1=-a_2=-b_2=q\sqrt{a}$ 
of Proposition~\ref{chugen2} where the series are specialized
to be very-well-poised. Then the first two terms on the right side of
\eqref{chugen2gl} vanish and we end up with the following ($r-1$)-term
transformation:
\begin{Proposition}[A very-well-poised $_{2r+2s}\psi_{2r+2s}$
Chu--Gasper--Karlsson--Minton type transformation]\label{chugen3}
Let $a$, $a_3,\dots,a_r$, $b_3,\dots,b_r$, and $h_1,\dots,h_s$
be indeterminate, let $m_1,\dots,m_s$
be nonnegative integers, $|m|=\sum_{i=1}^s m_i$, and suppose
that the series in \eqref{chugen3gl} are well-defined. Then
\begin{multline}\label{chugen3gl}
{}_{2r}\psi_{2r}\!\left[\begin{matrix}q\sqrt{a},-q\sqrt{a},b_3,\dots,b_{2r},
\{h_{\nu}\},\{\frac{aq^{1+m_{\nu}}}{h_{\nu}}\}\\
\sqrt{a},-\sqrt{a},\frac{a q}{b_3},\dots,\frac{a q}{b_{2r}},
\{\frac{aq}{h_{\nu}}\},\{h_{\nu}q^{-m_{\nu}}\}\end{matrix}\,;q,
\frac{a^{r-1}q^{r-2-|m|}}{b_3\dots b_{2r}}\right]\\
=\frac{(a,\frac q{a},a_4,\dots,a_r,\frac q{a_4},\dots,\frac q{a_r},
\frac{a_4}{a},\dots,\frac{a_r}{a},
\frac{aq}{a_4},\dots,\frac{aq}{a_r};q)_{\infty}}
{(\frac q{b_3},\dots,\frac q{b_{2r}},
\frac{a q}{b_3},\dots,\frac{a q}{b_{2r}},
\frac{a_4}{a_3},\dots,\frac{a_r}{a_3},\frac{a_3q}{a_4},\dots,
\frac{a_3q}{a_r};q)_{\infty}}\\\times
\frac{(\frac{a_3q}{b_3},\dots,\frac{a_3q}{b_{2r}},
\frac{a q}{a_3b_3},\dots,\frac{a q}{a_3b_{2r}};q)_{\infty}}
{(\frac{a_3a_4}{a},\dots,
\frac{a_3a_r}{a},\frac{aq}{a_3a_4},\dots,\frac{aq}{a_3a_r},
\frac{a_3^2q}{a},\frac{aq}{a_3^2};q)_{\infty}}
\prod_{i=1}^s\frac{(\frac{a_3q}{h_i},\frac{aq}{a_3h_i};q)_{m_i}}
{(\frac{aq}{h_i},\frac q{h_i};q)_{m_i}}\\\times
{}_{2r}\psi_{2r}\!\left[\begin{matrix}\frac{qa_3}{\sqrt{a}},
-\frac{qa_3}{\sqrt{a}},\frac{a_3b_3}{a},
\dots,\frac{a_3b_{2r}}{a},
\{\frac{a_3h_{\nu}}a\},\{\frac{a_3q^{1+m_{\nu}}}{h_{\nu}}\}\\
\frac{a_3}{\sqrt{a}},
-\frac{a_3}{\sqrt{a}},\frac{a_3q}{b_3},\dots,\frac{a_3q}{b_{2r}},
\{\frac{a_3q}{h_{\nu}}\},\{\frac{a_3h_{\nu}q^{-m_{\nu}}}a\}\end{matrix}\,;q,
\frac{a^{r-1}q^{r-2-|m|}}{b_3\dots b_{2r}}\right]\\
+\operatorname{idem}(a_3;a_4,\dots,a_r),
\end{multline}
where the series either terminate, or $|a^{r-1}q^{r-2-|m|}/b_1\dots b_{2r}|<1$,
for convergence.
\end{Proposition}

Finally, let us conclude this section stating an important
special case of Proposition~\ref{chugen3} explicitly. We take $r=3$,
$a_3\mapsto f$, $b_3\mapsto b$, $b_4\mapsto c$, $b_5\mapsto d$,
$b_6\mapsto e$, and obtain:
\begin{Corollary}\label{km3}
Let $a$ $b$, $c$, $d$, $e$, $f$, $h_1,\dots,h_s$
be indeterminate, let $m_1,\dots,m_s$
be nonnegative integers, $|m|=\sum_{i=1}^s m_i$, and suppose
that the series in \eqref{km3gl} are well-defined. Then
\begin{multline}\label{km3gl}
{}_{6+2s}\psi_{6+2s}\!\left[\begin{matrix}q\sqrt{a},-q\sqrt{a},b,c,d,e,
\{h_{\nu}\},\{\frac{aq^{1+m_{\nu}}}{h_{\nu}}\}\\
\sqrt{a},-\sqrt{a},\frac{a q}{b},\frac{a q}{c},\frac{a q}{d},
\frac{a q}{e},\{\frac{aq}{h_{\nu}}\},\{h_{\nu}q^{-m_{\nu}}\}\end{matrix}\,;q,
\frac{a^{2}q^{1-|m|}}{bcde}\right]\\
=\frac{(a,\frac q{a},\frac{fq}b,\frac{fq}c,\frac{fq}d,\frac{fq}e,
\frac{a q}{bf},\frac{a q}{cf},\frac{a q}{df},\frac{a q}{ef};q)_{\infty}}
{(\frac qb,\frac qc,\frac qd,\frac qe,\frac{a q}b,\frac{a q}c,
\frac{a q}d,\frac{a q}e,\frac{f^2q}{a},\frac{aq}{f^2};q)_{\infty}}
\prod_{i=1}^s\frac{(\frac{fq}{h_i},\frac{aq}{fh_i};q)_{m_i}}
{(\frac{aq}{h_i},\frac q{h_i};q)_{m_i}}\\\times
{}_{6+2s}\psi_{6+2s}\!\left[\begin{matrix}\frac{qf}{\sqrt{a}},
-\frac{qf}{\sqrt{a}},\frac{bf}{a},\frac{cf}{a},\frac{df}{a},\frac{ef}{a},
\{\frac{fh_{\nu}}a\},\{\frac{fq^{1+m_{\nu}}}{h_{\nu}}\}\\
\frac{f}{\sqrt{a}},-\frac{f}{\sqrt{a}},\frac{fq}{b},
\frac{fq}{c},\frac{fq}{d},\frac{fq}{e},
\{\frac{fq}{h_{\nu}}\},\{\frac{fh_{\nu}q^{-m_{\nu}}}a\}\end{matrix}\,;q,
\frac{a^{2}q^{1-|m|}}{bcde}\right],
\end{multline}
where the series either terminate, or $|a^{2}q^{1-|m|}/bcde|<1$,
for convergence.
\end{Corollary}

The special case $f\mapsto d$, $e\mapsto a/d$ of \eqref{km3gl}
is Chu's summation in \eqref{chugl}.

For multidimensional extensions of
Corollaries~\ref{km1}, \ref{km2} and \ref{km3} see \cite{schlkm},
where transformations of Chu--Gasper--Karlsson--Minton-type
are derived for multiple basic hypergeometric series
associated to the root system $A_n$.

\end{document}